\documentclass{amsart}
\usepackage{graphicx}
\usepackage{amssymb}
\usepackage{amscd}
\usepackage{latexsym}
\usepackage{amsfonts}
\usepackage{ragged2e}
\usepackage{amsmath}
\usepackage{float}
\usepackage[latin1]{inputenc}
\usepackage[english]{babel}

\newtheorem{proposition}{Proposition}%[section]
\newtheorem{lemma}{Lemma}%[section]
\newtheorem{theorem}{Theorem}%[section]
\newcommand{\eproof}{\begin{flushright} $\square$ \end{flushright}}

\let\eproof\endproof %better construction, the \eproof construction
\usepackage[draft]{fixme}

\newcommand\tildet[1]{r}

\newcommand{\Li}{\mathop{\fam0 Li}\nolimits}
\newcommand{\Log}{\mathop{\fam0 Log}\nolimits}

\newcommand{\ra}{\mathop{\fam0 \rightarrow}\nolimits}

\renewcommand{\epsilon}{\varepsilon}

\renewcommand{\phi}{\varphi}

\newcommand{\bbC}{\mathbb{C}}

\newcommand{\bbZ}{\mathbb{Z}}

\newcommand{\bbR}{\mathbb{R}}

\newcommand{\calB}{\mathcal{B}}

\newcommand{\calL}{\mathcal{L}}
\newcommand{\calM}{\mathcal{M}}

\newcommand{\calO}{\mathcal{O}}

\renewcommand{\Im}{\mathrm{Im}}
\renewcommand{\Re}{\mathrm{Re}}

\DeclareMathOperator{\Arg}{Arg}

\DeclareMathOperator{\Res}{Res}

\DeclareMathOperator{\U}{U}

\def\XXint#1#2#3{{\setbox0=\hbox{$#1{#2#3}{\int}$}
\vcenter{\hbox{$#2#3$}}\kern-.5\wd0}}

\begin{document}

\title{Resurgence Analysis of Meromorphic Transforms}

\author{J{\o}rgen Ellegaard Andersen}
\address{Center for Quantum Mathematics\\ 
  Danish Institute for Advanced Study\\
  University of Southern Denmark\\
  DK-5230, Denmark}

\email{jea-qm@mci.sdu.dk}

\thanks{This paper is partly a result of the ERC-SyG project, Recursive and Exact New Quantum Theory (ReNewQuantum) which received funding from the European Research Council (ERC) under the European Union's Horizon 2020 research and innovation programme under grant agreement No 810573.}
\maketitle

\begin{abstract}
We consider meromorphic transforms given by meromorphic kernels and study their asymptotic expansions under a certain rescaling. Under decay assumptions we establish the full asymptotic expansion in the rescaling parameter of these transforms and provide global estimates for error terms. We show that the resulting asymptotic series is Borel resummable and we provide formulae for the resulting resurgent function, which allows us to give formulae for the Stokes coefficients. A number of classical functions are obtained by applying such meromorphic transforms to elementary functions, of which, the Faddeev quantum dilogarithm, the Euler gamma function, the Riemann zeta function, the Gauss hypergeometric function and the Airy function are excellent examples of our general theory.
\end{abstract}

\section{Introduction}

In this paper, we study a certain type of meromorphic transforms which is very common in the study of classical functions. For a meromorphic function $h\in \calM(\bbC)$ we let $P_h$ denote its pole divisor. Further for $\gamma \in \bbC$ we define the rescaled meromorphic function $h_\gamma \in \calM(\bbC)$ by
$$ h_\gamma(z) = h(\gamma z).$$ 
Let $W\subset \bbC$ be an open connected subset and consider a meromorphic kernel $K \in \calM(W\times \bbC)$. We use the notation $ K_w = K|_{w \times \bbC}$ for $w\in  \bbC$, with the property that  $w \times \bbC$ is not contained in the pole divisor of $K$, which we will assume for all $w\in W$. 

Assume that we have a continuous family $\Gamma_w$,  $w\in W$ of curves, possibly non-compact and possibly with a fixed boundary, such that
$$ \Gamma_{w} \subset \bbC - P_{K_w}\ \ \forall w\in W.$$

We consider the following subspace $ \calM_{K,\Gamma}(\bbC) \subset \calM(\bbC)$ consisting  of the meromorphic functions $f\in \calM(\bbC)$ for which 
$$ \Gamma_{w} \subset \bbC - (P_{K_w} \cup P_{f_\gamma})$$
for all $(w, \gamma) \in  W\times U \subset \bbC\times \bbC$, where $U \subset \bbC$ is a non-empty connected open subset such  $0 \in \bar U - U$ and the integral
\begin{equation}\label{MO}
g_\gamma(w) = \int_{\Gamma_{w}} K(w,z) f_\gamma(z) dz
\end{equation}
exist for all $(w, \gamma) \in  W\times U$. This allows us to consider the following transforms
$$ A^\gamma_{K,\Gamma} : \calM_{K,\Gamma}(\bbC) \rightarrow \calO(W), \ \ \ A^\gamma_{K,\Gamma}(f) = g_\gamma.$$

We are in this paper interested in the asymptotic expansion of the transforms $A^\gamma_{K,\Gamma}$ as $\gamma$ tends to zero, e.g. the asymptotics of 
$g_\gamma$ as  $\gamma$ approached zero.

To this end we will assume that the following integrals  
$$ h_m(w) = \int_{\Gamma_{w}} K(w,z) z^m dz $$
exist for all\footnote{We might want to limit the pole order at $0$ of the meromorphic functions $f$ we allow, and thus only consider  our transform on the subspace
$ \calM_{K,\Gamma}^{n_0} = \{ f \in \calM_{K,\Gamma} \mid \text{ord}_0(f) \geq n_0\}$ and then we only need $m\geq n_0$.}  $m\in \bbZ$ and define $h_m\in \calO(W)$. We will further assume that the pole order of $K_w$ in zero is independent of $w$ and denote it $k_0 \in \bbZ_+\cap\{0\}$. 

Let now $n_0\in \bbZ_+\cup\{0\}$  be the pole order of $f$ in zero and $a_m \in \bbC$ the coefficients of the Laurent series of $f$ at zero
$$f(z) = \sum_{m=-n_0}^\infty a_m z^m$$
convergent for $z\in D(0,R_f)$, where $R_f$ is the minimal distance from $0\in \bbC$ to $P_f -\{0\}$. In particular, we can consider the smooth function $\kappa: (0,R_f) \rightarrow \bbR_+$ given by
\begin{equation}\label{kfn}
\kappa({r}) = \sum_{m=-n_0}^\infty |a_m| {r}^m, \ \ \ {r} \in (0,R_f),
 \end{equation}
 which we will use below in our analytic estimates.
We now define
$$ \tilde A^\gamma_{K,\Gamma}  :  \calM_{K,\Gamma}(\bbC) \rightarrow \calO(W)[\gamma^{-1},\gamma]]$$
given by
\begin{equation}\label{tildeA}
\tilde A^\gamma_{K,\Gamma}(f) = \sum_{m= - n_0}^\infty a_m h_m(w) \gamma^m.
\end{equation}
We introduce the following notation for the truncation at $n$ of the series (\ref{tildeA})
$$ 
\tilde A^{\gamma,n}_{K,\Gamma}(f) = \sum_{m= - n_0}^{n} a_m h_m(w) \gamma^m.
$$

In order to get analytic estimated on the reminder term in the asymptotic expansion of $A^\gamma_{K,\Gamma}$, we will make the following assumptions.

First of all, we will assume that $\Gamma_w$  can be continuously deformed (relative its boundary if non empty) in side $\bbC - (P_{K_w}\cup P_{f_\gamma} -\{0\})$ to a new contour $\tilde{\Gamma}_w$ consist of finitely many, say  $d$, smooth arc segments\footnote{If a smooth arc segment passes through zero, we split it at zero in two segments, simply to get the right definition of $d$ in relations to our estimates.}  $\tilde{\Gamma}_{w,j}$, $j=1, \ldots, d$ with the following property.

Let $\Theta(z)$ be the angle between the line through zero and $z$ and then the tangent line to $\tilde{\Gamma}_{w,j}$ at $z$, which is well defined along each of the smooth segments $\tilde{\Gamma}_{w,j}$. We then require that there exist a positive constant $b$ such that  in each smooth segment
\begin{equation}\label{cos}
|\cos(\Theta(z))| \geq b, \ \ \ \forall z\in \tilde{\Gamma}_{w,j}, \ j=1, \ldots d.
\end{equation}

We will also need decay conditions on $K$ and $f$.  In general we will proceed with the following assumptions.

There exist positive constants $c, \delta_1, \delta_2$ such that
$$ \delta := \delta_1 - \delta_2 >0$$
together with $ c_w \in \bbR_+$ parametrized by $w\in W$ and $\tilde{c}_\gamma \in [c,\infty)$ parametrized by $\gamma\in U$,  such that 
\begin{equation} \label{Mainestimateinfty}
\mid K(w,z) \mid \leq c_w e^{-\delta_1 | z |} |z|^{-k_0}  \ \ \ \forall (w,z) \in W \times (\tilde{\Gamma}_w - \{0\})
\end{equation}
\begin{equation} \label{Mainestimateinftyf}
|f_\gamma(z)| \leq \tilde{c}_\gamma e^{\delta_2 |z| }|\gamma z|^{-n_0}  \ \ \ \forall (\gamma, z) \in U\times (\tilde{\Gamma}_w - \{0\}).
\end{equation}

The specifics of these norm estimated together with the assumptions on $\tilde{\Gamma}_w$ are tailored to establish the global control for all $(w,\gamma) \in W\times U$ in (\ref{MEP}) in the following theorem.

Let  
$$ C_{k_0-1} = \frac{d |\gamma| }{bc\delta}c'_{k_0}, \ \ \ C_n = \frac{2d }{bc\sqrt{\pi(1+ 2(n-k_0))}}c'_n, \ \ \ n\geq k_0$$
where
$$ c'_n = \inf_{0<{r} < R_f} \frac{1}{{r}^{n}} \left(cr^{-n_0} +  \kappa(r)\right), \ \ \ n\geq k_0.$$
\begin{theorem}\label{Thm1}
Under the above assumptions, we have an asymptotic expansion in the Poincare sense to all orders
\begin{equation}\label{gamma_expansion}
A^\gamma_{K,\Gamma}(f) \sim \tilde A^\gamma_{K,\Gamma}(f).
\end{equation}
In fact we have the following estimates for  all $n \geq k_0-1$ 
\begin{equation}\label{MEP}
\mid A^\gamma_{K,\Gamma}(f)(w)  -  \tilde A^{\gamma,n}_{K,\Gamma}(f)(w) \mid \leq C_n c_w\tilde{c}_\gamma|\gamma|^{n} \delta^{-(n-k_0+1)}  (n-k_0+1)! 
\end{equation}
for all $(w,\gamma) \in W\times U$.
\end{theorem}

This theorem together with the Theorems \ref{Thm1.5} and \ref{Thm1.75} in Section \ref{AE}, which covers the case $-n_0\leq n< k_0$ under a slight modification of the assumptions (\ref{cos}), (\ref{Mainestimateinfty}) and (\ref{Mainestimateinftyf}), are proved in Section \ref{AE}. The basic principle is to spilt the integral (\ref{MO}) into the sum of the integral over $\Gamma^{\gamma,0}_w(\tilde{r}) = \Gamma^{\gamma,0}_w\cap \overline{D(0,\tilde{r})}$ and over $\Gamma^{\gamma,\infty}_w(\tilde{r}) = \Gamma^{\gamma}_w \cap (\bbC- D(0,\tilde{r}))$ for some positive $\tilde{r}$. The first of these integrals is simply estimated using (\ref{Mainestimateinfty}) and (\ref{Mainestimateinftyf}). For the part of the integral over $\Gamma^{\gamma,0}_w(\tilde{r})$, we require that ${r}=|\gamma|\tilde{r} < R_f$ and then we can use the Laurent expansion for $f_\gamma$ to get the needed estimate on that part of the integral to then establish (\ref{MEP}). The estimate  (\ref{MENs}) in Theorem \ref{Thm1.5} is established along similar lines under slightly different assumptions and the same is the case for Theorem \ref{Thm1.75}, where we limit the condition (\ref{cos}) to the part of $\tilde{\Gamma}_w$ outside $D(0,\tilde{r})$, to obtain yet another estimate. All these different assumptions are handy when dealing with various different examples.

If we have that $R_f = \infty$, then we can make milder conditions dealing with a much simpler case whose discussion we defer to Section \ref{Mild}.

We are further interested in the resurgence properties of the Borel resummation of the resulting possibly divergent series $\tilde{A}^\gamma_{K,\Gamma}(f)$ in $\gamma$. 
 To this end we introduce  $g^-_\gamma \in\calM(W)[\gamma]$, given by
$$g^-_{\gamma}(w) = \sum_{m= -n_0}^{k_0-1} a_{m }h_{m}(w) \gamma^{m},$$
and $g_\gamma^+ \in \gamma \calM(W)[[\gamma]]$ defined by
$$g^+_{\gamma}(w) =  \sum_{m= k_0}^\infty a_m h_m(w) \gamma^m.$$

We will typically expect that (\ref{gamma_expansion}) is only an asymptotic expansion, e.g.  that the series $g^+_{\gamma}(w)$ is divergent, as we will see in various examples in Section \ref{Faddeev} and \ref{Examples}.
We therefore consider the formal Borel transform 
$$ \calB : \gamma\calO(W)[[\gamma]] \rightarrow \calO(W)[[\xi]]$$
determined by
$$ \calB(\gamma^{m})  = \frac{\xi^{m-1}}{(m-1)!}$$
and formally extended linearly over $\calO(W)$. We recall that $\calB$ applied to a divergent Gevrey-1 power series gives a power series with positive radius of convergence.

Let $\bbR_+^\theta = e^{i\theta} \bbR_+$ for $\theta \in \bbR$ and consider the Laplace transform
$$ \calL_\theta(\psi)(\gamma) = \int_{\bbR_+^\theta} e^{-\xi/\gamma} \psi(\xi) d\xi,$$
which is certainly well defined for measurable $\psi$ defined on $\bbR_+^\theta$ provided there exist a positive constant $C$, real $\alpha$ and a non-negative integer $m$ such that
$$ |\psi(\xi) | \leq C e^{\alpha |\xi|} |\xi|^m, \forall \xi\in \bbR_+^\theta$$
and
$$
\frac{1}{|\gamma|} \cos(\theta - \theta_\gamma) >\alpha 
$$
where
$$
\theta_\gamma := \text{Arg}(\theta).$$
Let now 
$$U_{\theta,\alpha}= \left\{ \gamma \in U\bigg\vert \cos(\theta - \theta_\gamma) > |\gamma| \alpha \right\} \ \ \text{ and } \ \ \U_{\theta}= \left\{ \gamma \in U\bigg\vert \cos(\theta - \theta_\gamma) > 0\right\}.$$
We observe that $\calL_\theta (\xi^{m-1})(\gamma)$ , $m\in \bbZ_+$, is well defined for $\gamma \in  U_{\theta,\alpha}$ for all $\alpha \in \bbR_+\cup \{0\}$, thus
$$\calL_\theta (\xi^{m-1}) \in \calO( U_{\theta})$$
and that 
$$\calL_\theta  \circ \calB (\gamma^m) = \gamma^m \ \  \forall m\in \bbZ_+.$$

We will now state our assumptions which will guarantee that $\calB(g^+_{\gamma}(w))$ is convergent for all $w\in W$.

To that end introduce $\phi \in \calM(\bbC)$ given by
$$ \phi(z) = f(z) - \sum_{m=-n_0}^{k_0-1} a_{m} z^{m}$$
and $\phi_z \in \calM(\bbC)$ given by
$$ \phi_z(\gamma) = \phi(\gamma z).$$

We will now assume there exist a continuous deformation of $\tilde{\Gamma}_w$ in side $\bbC - P_{K_w}$ to a new contour
$$ \overline{\Gamma}_w \subset \bbC - (P_{K_w}\cup \{0\}), \ \ \ \forall w \in W$$
such that there exist $r(w)\in  \bbR_+$ with the property that (\ref{cos}), (\ref{Mainestimateinfty}) and  (\ref{Mainestimateinftyf})  holds along 
\begin{equation} \label{Asump3} 
\overline{\Gamma}_w^0(r(w)) = \overline{\Gamma}_w \cap D(0, r(w)).
\end{equation}

\begin{theorem}\label{Thm2}
Fix $\theta \in \bbR$. Suppose there exist
$$ \alpha_\theta: W \rightarrow \bbR$$
and finitely many constant $C_i$ and non-negative integers $m_i$ and  $V_{w,\theta}\subset \bbC$ an open subset containing the half line $\overline{\bbR_+^\theta}$ with the property that for all $w \in W$
\begin{equation}\label{hol}
\calL_\theta^{-1}(\phi_z) \in \mathcal{O}(V_{w,\theta}) \ \ \ \forall z\in \overline{\Gamma}_w
\end{equation}
and
\begin{equation}\label{intest}
| \calL_\theta^{-1}(\phi_z)(\xi) | \leq e^{\alpha_\theta(w) |\xi| } \sum_{i} C_i |\xi|^{m_i} \ \ \ \forall \xi \in V_{w,\theta}.
\end{equation}
Further we assume
\begin{equation}\label{alphalb}
\gamma \in U_{\theta,\alpha_\theta(w)} \ \ \ \forall w\in W.
\end{equation}
Then we have that the formal series $g^+_\gamma$ is Borel summable and
\begin{equation}\label{Borel-resum}
 B_w(\xi) = \calB(g^+_\gamma(w))(\xi)
\end{equation}
 is well defined for $\xi \in V_{w,\theta}$ and  $B_w \in \calO(V_{w,\theta})$ is given by 
\begin{equation}\label{resurgence-Borel}
 B_w(\xi) = \int_{\overline{\Gamma}_{w}} K(w,z) \calL_\theta^{-1}(\phi_z)(\xi) dz.
\end{equation}
Further we have that 
\begin{equation}\label{Laplace-Borel}
g_\gamma (w) =  g^-_\gamma(w) + \calL_\theta(B_w)(\gamma) \ \ \ \forall \gamma \in U_{\theta,\alpha_\theta(w)} \ \ \forall w\in W.
\end{equation}
\end{theorem}

The proof of this theorem is provided in Section \ref{Borel}. We first establish the formula  (\ref{Laplace-Borel}) under the assumptions stated in the theorem and then we can derive the rest by properties of the Laplace transform.

The resurgence properties and Stokes phenomenon of $B_w$ can now be studied. In general we see that $\calL_\theta(B_w)$ will be constant on sectors where $B_w$ has no poles and jump when $\bbR_+^\theta$ hits poles of  $B_w$ along a set of directions  
$\theta_{J_w} = \{ \theta_j \mid j \in J_w\}$ index by a set $J_w$, $w\in W$.
Let the jump at $\theta_j$ be denote $\Delta_{\theta_j} (\calL_\theta(B_w))$. 

\begin{theorem}\label{Thm3}
Let $w\in W$. Assume there exist $\epsilon >0$ such that $\calL_\theta(\calB_w)$ is well defined for $\theta \in (\theta_j, \theta_j \pm \epsilon]$. Assume further the poles of $B_w$ in this sector $\theta \in [\theta_j-\epsilon, \theta_j + \epsilon]$ are on the line $\bbR_+^\theta$ and that there exist a sequence $R_n \in \bbR_+$, which converges to infinity as $n \rightarrow \infty$ such that
\begin{equation} \label{decayc}
 \lim_{n\rightarrow \infty}   R_n \sup_{\xi \in \Gamma_n} \vert e^{-\xi/\gamma}B_w(\xi)) \vert = 0,
 \end{equation}
where
$$ \Gamma_n = \left\{ R_n e^{i\theta} \mid \theta \in [\theta_j-\epsilon, \theta_j + \epsilon] \right\}$$
and $R_n e^{i\theta_j}$ is not a pole of $B_w$.
Then
$$\Delta_{\theta_j} (\calL_\theta(B_w)) = 2 \pi i\sum_{p\in P_{B_w} \cap \bbR_{\theta_j}^+} \Res_{\xi = p} (e^{-\xi/\gamma} B_w(\xi)). $$
\end{theorem}

Let now the principal part of $\phi_z$ at $p\in P_{\phi}$ be denoted 
$$\tilde{\phi}_{z,p}(\gamma) = \sum_{m=1}^{n_\varphi} \frac{b_{p,m}}{(\gamma z - p)^m},$$
where we assume that the pole order is universally bounded by some integer $n_\varphi $ independent of $z \in \tilde{\Gamma}_w$ and $w\in W$.
We observe that if $p\neq 0$ and we let
$$ \tilde{\phi}'_{z,p}(\gamma) = -\sum_{m=1}^{n_\varphi} \frac{b_{p,m}}{p^m}  \sum_{l=0}^{m-1}  {{m}\choose{l}} (-1)^l\sum_{l'=0}^l  {{l}\choose{l'}} \left(\frac{z}{p}\right)^{m-l+l'} \frac{1}{(\frac{1}{\gamma} - \frac{z}{p})^{m-l+l'}}$$
then
$$\tilde{\phi}_{z,p}(\gamma) = \sum_{m=1}^{n_\varphi} \frac{b_{p,m}}{p^m} (-1)^m   +  \tilde{\phi}'_{z,p}(\gamma).$$
We observe that the series 
$$ \tilde{\phi}'_z(\gamma) = \sum_{p\in P_{\phi}} \tilde{\phi}'_{z,p}(\gamma)$$
is convergent uniformly on 
$$O_{\epsilon,z} = \left\{ \gamma \in \bbC \ \bigg\vert  \ \bigg\vert \frac{1}{\gamma} - \frac{z}{p}\bigg\vert > \epsilon \ \ \forall p\in P_{\phi} \right\}$$
for all $\epsilon >0$, when we partially ordering the points in $P_{\phi}$ by increasing distance to $0$ provided that
\begin{equation}\label{Asump1}
\sum_{p\in P_\phi \cap (D(0,n)-D(0,n-1))} \vert b_{p,m}\vert \leq \tilde{C} |p|^c \ \ \ \text{ for } \ m=1, \ldots, n_\phi
\end{equation}
where $c <1$ and $\tilde{C}$ is some constant and $0 \notin P_\phi$.
We will further assume that there exist an entire function $\psi_z$ parametrized by  $z \in \overline{\Gamma}_w$ for $w\in W$ such that
\begin{equation}\label{Asump2}
\phi_z(\gamma) = \tilde{\phi}'_z(\gamma) + \psi_z(\gamma) \ \ \ \forall \gamma \in O_{\epsilon,z} \  \forall z \in \overline{\Gamma}_w \ \forall w\in W.
\end{equation}
We will see in examples that one actually sometimes get that $\psi_z = 0$.
Introduce the following notations for the arguments
$$
\theta_\gamma := \text{Arg}(\gamma) , \ \  \theta_z := \text{Arg}(z) \ \ \text{ and } \ \  \theta_p : = \text{Arg}(p).
$$
We will further need the following conditions
\begin{equation}\label{cosrel}
\cos(\theta - \theta_\gamma) >0 \ \ \ \text{ and } \ \ \
 \gamma \in U_{\theta,\alpha_\theta(w)}
\end{equation}
where 
$$  \alpha_\theta(w) = \frac{1}{r_m} \sup_{z\in \overline{\Gamma}_w} |z|\cos(\theta + \theta_z -\theta_p). $$
for all  $w\in W$ and
\begin{equation}\label{not0}
 r_m = \min_{p\in P_\phi} |p| >0.
 \end{equation}
 Note that (\ref{not0}) is equivalent to $0\notin P_\phi$.

As we will see in examples, the conditions (\ref{Asump1}), (\ref{Asump2}), (\ref{cosrel}) and (\ref{not0}) are geometric conditions on $P_\phi$, which are easy to check.

\begin{theorem}\label{Thm5} 
Fix $\theta \in \bbR$. Suppose there exist
$$ \alpha'_\theta : W \rightarrow \bbR$$
and finitely many constants $C_i$ and nonnegative integers $m_i$ and further $V_{w,\theta}\subset \bbC$ an open subset containing the half line $\overline{\bbR_+^\theta}$ with the property that for all $w \in W$
\begin{equation}\label{phol}
\calL_\theta^{-1}(\psi_z) \in \mathcal{O}(V_{w,\theta}) \ \ \ \forall z\in \overline{\Gamma}_w
\end{equation}
and
\begin{equation}\label{pintest}
 | \calL_\theta^{-1}(\psi_z)(\xi) |  \leq e^{\alpha'_\theta(w) |\xi| } \sum_i C_i |\xi|^{m_i} \ \ \ \forall \xi \in V_{w,\theta}.
\end{equation}
Further assume (\ref{Asump1}), (\ref{Asump2}), (\ref{cosrel}) and (\ref{not0}).
Then we conclude that 
\begin{eqnarray}\label{B_w-sum}
B_w(\xi) &= & -\sum_{p\in P_{\phi}}  \sum_{m=1}^{n_\varphi} \frac{b_{p,m}}{p^m}  \nonumber \sum_{l=0}^{m-1}  {{m}\choose{l}}(-1)^l \sum_{l'=0}^l  {{l}\choose{l'}}  \\
& & \phantom{jjjjj} \left( \int_{\overline{\Gamma}_{w}} K(w,z) e^{\frac{z\xi}{p}} \left(\frac{z}{p}\right)^{m-l+l'}  dz \right) \frac{\xi^{m-l+l'-1}}{(m-l+l'-1)!}  \nonumber\\
& & + \int_{\overline{\Gamma}_{w}} K(w,z) \calL_\theta^{-1}(\psi_z)(\xi) dz.
\end{eqnarray}
for all $w \in W$ and $\xi \in V_{w,\theta}$.
\end{theorem}

\noindent The proofs of Theorem \ref{Thm3} and \ref{Thm5} are provided in Section \ref{Stokes}. The formula (\ref{B_w-sum}) is ideally suited to understanding the resurgence properties of $B_w$ and the Stokes coefficients, since it is in examples easy to determined what maximal open subset $B_w$ extends to from the expression (\ref{B_w-sum}), as we will in particular see in the case of Faddeev quantum dilogarithm in Section \ref{Faddeev}.

In the final two sections of this paper we provide several examples of our main theorems stated above, which in particular include the case of the Faddeev quantum dilogarithm, the Euler gamma function, the Riemann zeta function, the Hurwitz zeta function, the Gauss hypergeometric function and the Airy function. 

\textbf{Acknowledgements.} The author thanks  Stavros Garoufalidis and Rinat Kashaev  for helpful discussion. In fact, some of the present work, in particular the applications of our main theorems to Faddeev's quantum dilogarithm was produced by the author of this paper about five years ago in preparation for a joint project with Garoufalidis and Kashaev concerned with the asymptotics of our Quantum Teichm\"{u}ller theory constructed by the author of this paper jointly with Kashaev \cite{AK1,AK2,AK3,AK4,AK5,AK6, AN,AM,AMa} and the KLV invariant presented in \cite{KLV}. The author is partially supported by the ERC synergy grant "ReNewQuantum".

\section{Proof of the asymptotic expansion}

\label{AE}

To analyse the asymptotics of the transform $A^\gamma_{K,\Gamma}$ we introduce the function $\varphi_n\in \calM(\bbC)$ given by
$$\varphi_n(z) = f(z) - \sum_{m=-n_0}^{n} a_m z^m.$$
When we assume that $n\geq k_0-1$, we get that 
$$z\mapsto K(w,z) \varphi_n(\gamma z) \text{ is regular } \forall z\in D(0,\min(R_K(w),\frac{R_f}{|\gamma|}))$$
for all $w\in W$, where $R_K(w)$ is the minimal distance from zero to $P_{K_w} - \{0\}$. Thus we see that
$$ A^\gamma_{K,\Gamma}(f)(w)  -  \tilde A^{\gamma,n}_{K,\Gamma}(f)(w) = \int_{\tilde{\Gamma}_w} K(w,z) \varphi_n(\gamma z)dz$$
and therefore
$$  \mid A^\gamma_{K,\Gamma}(f)  -  \tilde A^{\gamma,n}_{K,\Gamma}(f) \mid \leq \int_{\tilde{\Gamma}_w} |K(w,z)| |\varphi_n(\gamma z)| |dz|$$
Let now $\tilde{r} <\frac{R_f}{|\gamma|}$ and for each $j=1,\ldots d$ consider the decomposition
$$ \tilde{\Gamma}_{w,j} = \tilde{\Gamma}^{0}_{w,j}(\tilde{r}) \cup \tilde{\Gamma}^{\infty}_{w,j}(\tilde{r}),$$
where 
$$ \tilde{\Gamma}^{0}_{w,j}(\tilde{r}) = \tilde{\Gamma}_{w,j} \cap \overline{D(0,\tilde{r})} \ \  \text{ and } \ \ \tilde{\Gamma}^{\infty}_{w,j}(\tilde{r}) = \tilde{\Gamma}_{w,j} \cap (\bbC -D(0,\tilde{r}))$$
and we let 
$$  \tilde{\Gamma}^{0}_{w}(\tilde{r})  =  \bigcup_{j=1}^d \tilde{\Gamma}^{0}_{w,j}(\tilde{r}), \ \ \ \tilde{\Gamma}^{\infty}_{w}(\tilde{r})  =  \bigcup_{j=1}^d \tilde{\Gamma}^{\infty}_{w,j}(\tilde{r}).$$
We immediately get the estimate
\begin{eqnarray*} 
\mid A^\gamma_{K,\Gamma}(f)(w)  -  \tilde A^{\gamma,n}_{K,\Gamma}(f)(w) \mid & \leq & \int_{\tilde{\Gamma}_w^{\infty}(\tilde{r})} |K(w,z)| |f(\gamma z)| |dz|\\
& + & \sum_{m=-n_0}^{n} |a_m| \int_{\tilde{\Gamma}_w^{\infty}(\tilde{r})} |K(w,z)| |\gamma z|^m |dz|\\
& + & \int_{\tilde{\Gamma}_w^{0}(\tilde{r})} |K(w,z)| |\varphi_n(\gamma z)| |dz|.\\
\end{eqnarray*} 
We now provide estimates for each of these three summands. First of all, it follows directly from our main estimates (\ref{Mainestimateinfty}) and (\ref{Mainestimateinftyf}) that
$$\int_{\tilde{\Gamma}_w^{\infty}(\tilde{r})} |K(w,z)| |f(\gamma z)| |dz| \leq c_w \tilde{c}_\gamma |\gamma|^{-n_0}  \int_{\tilde{\Gamma}_w^{\infty}(\tilde{r})}  e^{-\delta |z|} |z|^{-k_0-n_0} |dz|$$
and 
$$\int_{\tilde{\Gamma}_w^{\infty}(\tilde{r})} |K(w,z)| |\gamma z|^m |dz| \leq  c_w |\gamma|^m\int_{\tilde{\Gamma}_w^{\infty}(\tilde{r})}  e^{-\delta_1 |z|} |z|^{m-k_0} |dz|$$
For fixed $w$, we  now let $ \gamma_j: [t_{j,0}, t_{j,1}] \ra \tilde{\Gamma}^{\infty}_{w,j}(\tilde{r})$, $ j=1,\ldots d$, be a smooth parametrization of $\tilde{\Gamma}^{\infty}_{w,j}(\tilde{r})$ and for those $j$ such that $ \tilde{\Gamma}^{\infty}_{w,j}(\tilde{r}) = \emptyset$ we just set $t_{j,0} = t_{j,1}$. For some $j$ we might have $t_{j,1} = \infty$ and we can assume that $|\gamma_j(t_{j,0})| \leq |\gamma_j(t_{j,1})|$. If $\tilde{\delta} >0$, we see that
$$
\int_{\tilde{\Gamma}_w^{\infty}(\tilde{r})}  e^{-\tilde{\delta} |z|} |z|^{l} |dz| \leq \sum_{j=1}^{d}\int_{t_{j,0}}^{t_{j,1}} e^{-\tilde{\delta} |\gamma_j(t)|} |\gamma_j(t)|^l |\gamma'_j(t)| dt.$$
But since we have the lower bound $b$ on $|\cos(\Theta(z))|$ from (\ref{cos}) for $z\in \tilde{\Gamma}_{w,j}^{\infty}(\tilde{r})$, we see that
$$\int_{t_{j,0}}^{t_{j,1}} e^{-{\tilde{\delta}} |\gamma_j(t)|} |\gamma_j(t)|^l |\gamma'_j(t)| |dt| \leq \frac{1}{b} \int_{|\gamma(t_{j,0})|}^{|\gamma(t_{j,1})|} e^{-{\tilde{\delta}} s} s^l ds,$$
since we can make the substitution $s= |\gamma_j(t)|$ and then use that
$$ |\gamma_j(t)|' = \frac{\Re(\gamma_j(t)) \Re(\gamma'_j(t)) + \Im(\gamma_j(t))\Im(\gamma'_j(t))}{|\gamma_j(t)|} = \frac{|\gamma_j(t)||\gamma'_j(t)| \cos(\Theta(\gamma_j(t)))}{|\gamma_j(t)|},$$
which in absolute value is bounded below by $|\gamma'_j(t)|b$.

For $\tilde{r}_1\geq \tilde{r}_0 \geq \tilde{r}$ and $l\in \bbZ_-\cup \{0\}$ we have that
$$\int_{\tilde{r}_0}^{\tilde{r}_1} e^{-{\tilde{\delta}} s} s^l ds \leq \frac{1}{{\tilde{\delta}}}e^{-{\tilde{\delta}} \tilde{r}} \tilde{r}^{l}$$
and for $l\in \bbZ_+$ we have that
\begin{eqnarray*}
\int_{\tilde{r}_0}^{\tilde{r}_1} e^{-{\tilde{\delta}} s} s^l ds & \leq & \int_{\tilde{r}}^{\tilde{r}_1} e^{-{\tilde{\delta}} s} s^l ds\\
& = & \sum_{l'=0}^{l}  {\tilde{\delta}}^{-(l'+1)} \frac{l!}{(l-l')!} (e^{-{\tilde{\delta}} \tilde{r}} \tilde{r}^{l-l'} - e^{-{\tilde{\delta}} \tilde{r}_1} \tilde{r}_1^{l-l'})\\
& \leq &\sum_{l'=0}^{l}  {\tilde{\delta}}^{-(l'+1)} \frac{l!}{(l-l')!} e^{-{\tilde{\delta}} \tilde{r}} \tilde{r}^{l-l'}.
\end{eqnarray*}
Thus we have proved the following two lemmas, which provides the estimates we need in regards to the $\tilde{\Gamma}^{\infty}_w(\tilde{r})$-part of the integral.

\begin{lemma}\label{inftyf}
We have that
$$\int_{\tilde{\Gamma}_w^{\infty}(\tilde{r})} |K(w,z)| |f(\gamma z)| |dz| \leq \frac{d \tilde{c}_\gamma  c_w }{b} \delta^{-1} |\gamma|^{-n_0}e^{-\delta \tilde{r}} \tilde{r}^{-(k_0+n_0)}$$
\end{lemma}
For the second integral over $\tilde{\Gamma}_w^{\infty}(\tilde{r})$, we obtain the following estimates. 
\begin{lemma}\label{inftyznk0}
Provided that $n\geq k_0-1$, it follows that
\begin{eqnarray*}
\sum_{m=-n_0}^{n} |a_m| \int_{\tilde{\Gamma}_w^{\infty}(\tilde{r})} |K(w,z)| |\gamma z|^m |dz| & \leq & \frac{d c_w }{b}\delta_1^{-1} e^{-\delta_1 \tilde{r}}\sum_{m=-n_0}^{k_0-1} |a_{m}| |\gamma|^{m}  \tilde{r}^{m-k_0}\\
& +& \frac{d c_w }{b}e^{-\delta_1 \tilde{r}}\sum_{m=k_0}^{n} |a_m| |\gamma|^m\\
& & \phantom{hhhhh} \sum_{l=0}^{m-k_0}  \delta_1^{-(l+1)} \frac{(m-k_0)!}{(m-k_0-l)!}  \tilde{r}^{m-k_0-l}.
\end{eqnarray*}
\end{lemma}

Let us now attend to the $\tilde{\Gamma}^{0}_w(\tilde{r})$-part of the integral. We introduce the notation
$$ \kappa_n(r) = \sum_{m=n+1}^\infty |a_m| {r}^m.$$
\begin{lemma}\label{0phiGk_0}
Provided that ${r} = \tilde{r}|\gamma|$ satisfies ${r}<R_f$, we have for $n \geq k_0$ 
\begin{equation*}%\label{Near0EP} 
\int_{\tilde{\Gamma}_w^{0}(\tilde{r})} |K(w,z)| |\varphi_n(\gamma z)| |dz| \leq   \frac{dc_w  \kappa_{n+1}({r}){r}^{-n} }{b\sqrt{2\pi(n-k_0+1)}}|\gamma|^{n}  \delta_1^{k_0 -n -1} (n-k_0+1)!
\end{equation*}
and for $n = k_0-1$ 
\begin{equation*}%\label{Near0EP} 
\int_{\tilde{\Gamma}_w^{0}(\tilde{r})} |K(w,z)| |\varphi_n(\gamma z)| |dz| \leq   \frac{dc_w  \kappa_{k_0}({r}){r}^{-n-1} }{b\delta_1}|\gamma|^{n+1} 
\end{equation*}
\end{lemma}
\proof
We introduce the meromorphic function $\tilde{\varphi}_n \in \calM(\bbC)$ given by
$$ \tilde{\varphi}_n(\tilde{z}) = \tilde{z}^{-n-1} \varphi_n(\tilde{z}).$$
We get that
$$ |\tilde{\varphi}_n(\tilde{z})| \leq {r}^{-n-1} \sum_{m=n+1}^\infty |a_m| {r}^m = {r}^{-n-1} \kappa_{n+1}({r})$$
for all $\tilde{z}\in D(0,{r})$.
From this we conclude that
\begin{eqnarray*} 
\int_{\tilde{\Gamma}_w^{0}(\tilde{r})} |K(w,z)| |\varphi_n(\gamma z)| |dz| & = &  \int_{\tilde{\Gamma}_w^{0}(\tilde{r})} |K(w,z)| |\tilde{\varphi}_n(\gamma z)| |\gamma z|^{n+1} |dz|\\
& \leq &  c_w  \kappa_{n+1}({r}){r}^{-n-1} |\gamma|^{n+1}  \int_{\tilde{\Gamma}_w^{0}(\tilde{r})} e^{-\delta_1 |z|} |z|^{n-k_0+1} |dz|.
\end{eqnarray*}
since $z\in \tilde{\Gamma}_w^{0}(\tilde{r})$ implies that  $|z| \leq \tilde{r}$ and so $|\gamma z| \leq {r}$.
Thus we obtain
$$\int_{\tilde{\Gamma}_w^{0}(\tilde{r})} |K(w,z)| |\varphi_n(\gamma z)| |dz| \leq \frac{d c_w \kappa_{n+1}({r}){r}^{ -n-1}}{b} |\gamma|^{n+1}  \int_{0}^{\tilde{r}} e^{-\delta_1 s} s^{n-k_0+1} ds$$
and so for $n= k_0-1$ we get that
$$
\int_{0}^{\tilde{r}} e^{-\delta_1 s} s^{n-k_0+1} ds  \leq   \delta_1^{-1} 
$$
which gives the stated result in case $n=k_0-1$. In the case $n\geq k_0$ we use Robbin's lower bound for $m!$
 $$ \sqrt{2\pi m} \sup_{x\in \bbR_+} x^m e^{-\delta_1 x} = \sqrt{2\pi m} \delta_1^{-m} m^m e^{-m}\leq \delta_1^{-m} m!$$
 to get that
 $$
\int_{0}^{\tilde{r}} e^{-\delta_1 s} s^{n-k_0+1} \leq   \frac{\tilde{r} }{\sqrt{2\pi(n-k_0+1)}}\delta_1^{k_0 -n -1} (n-k_0+1)!
 $$
 which complets the proof.
\eproof

We now completing the proof of Theorem \ref{Thm1}  based on Lemma \ref{inftyf}, \ref{inftyznk0} and \ref{0phiGk_0}.

\proof 
We assume that $n \geq k_0-1$ and further that ${r} = \tilde{r}|\gamma| < R_f$. With reference to Lemma \ref{inftyf}, we compute 
\begin{eqnarray*} 
\delta^{-1} |\gamma|^{-n_0}e^{-\delta r} \tilde{r}^{-(k_0+n_0)} & \leq & r^{-n_0} |\gamma|^{ n} \delta^{-(n-k_0 +1)} (n-k_0+1)!\\
 & & \frac{{r}^{-n}}{(n-k_0+1)!}\left(\delta \frac{{r}}{|\gamma|}\right)^{n-k_0} e^{-\delta \frac{{r}}{|\gamma|}}
  \end{eqnarray*}
Thus for $n=k_0-1$ we get that
$$\delta^{-1} |\gamma|^{-n_0}e^{-\delta r} \tilde{r}^{-(k_0+n_0)}  \leq  {r}^{-n_0}  r^{-n-1} |\gamma|^{ n+1} \delta^{-1},$$
for $n=k_0$
$$ \delta^{-1} |\gamma|^{-n_0}e^{-\delta r} \tilde{r}^{-(k_0+n_0)}  \leq r^{-(n_0+n)} |\gamma|^{ n} \delta^{-(n-k_0 +1)} (n-k_0+1)!$$
and for $n>k_0$
$$\delta^{-1} |\gamma|^{-n_0}e^{-\delta r} \tilde{r}^{-(k_0+n_0)}  \leq \frac{{r}^{-(n+n_0)}}{\sqrt{2\pi(n-k_0+1)}} |\gamma|^{ n} \delta^{-(n-k_0 +1)} (n-k_0+1)!$$
 using again Robbin's lower bound for $m!$.
We now introduce the notation
$$\kappa_{n,n'}(r) = \sum_{m=n}^{n'} |a_{m}| {r}^{m.}$$
 With reference to Lemma \ref{inftyznk0}, using $\delta_1 > \delta$, we compute 
\begin{eqnarray*} 
 \delta^{-1} e^{-\delta \tilde{r}}\sum_{m=-n_0}^{k_0-1} |a_{m}| |\gamma|^{m}  \tilde{r}^{m-k_0} & \leq & \delta^{-1} e^{-\delta \frac{{r}}{|\gamma|}} |\gamma|^{k_0}{r}^{-k_0} \sum_{m=-n_0}^{k_0-1} |a_{m}|  {r}^{m}\\
 & \leq & {r}^{-k_0} \kappa_{-n_0,k_0-1}({r})|\gamma|^{n}\delta^{-(n-k_0 +1)} (n-k_0+1)! \\
 & & \frac{\left(\delta\frac{{r}}{|\gamma|}\right)^{n-k_0}e^{-\delta \frac{{r}}{|\gamma|}}}{(n-k_0+1)!} {r}^{-(n- k_0)} 
 \end{eqnarray*}
 which for $n=k_0-1$
 $$ \delta^{-1} e^{-\delta \tilde{r}}\sum_{m=-n_0}^{k_0-1} |a_{m}| |\gamma|^{m}  \tilde{r}^{m-k_0}  \leq  {r}^{-k_0}\kappa_{-n_0,k_0-1}({r})|\gamma|^{n+1}\delta^{-1},$$
 and for $n=k_0$
  $$ \delta^{-1} e^{-\delta \tilde{r}}\sum_{m=-n_0}^{k_0-1} |a_{m}| |\gamma|^{m}  \tilde{r}^{m-k_0}  \leq  {r}^{-n}\kappa_{-n_0,k_0-1}({r})|\gamma|^{n}\delta^{-(n-k_0+1)} (n-k_0+1)!$$
 and for $n>k_0-1$
$$ \delta^{-1} e^{-\delta \tilde{r}}\sum_{m=-n_0}^{k_0-1} |a_{m}| |\gamma|^{m}  \tilde{r}^{m-k_0}  \leq  \frac{{r}^{-n}\kappa_{-n_0,k_0-1}({r})}{\sqrt{2\pi(n-k_0)}(n-k_0+1)}|\gamma|^{n}\delta^{-(n-k_0 +1)} (n-k_0+1)! $$
and further if $n\geq k_0$
 \begin{eqnarray*} 
\mbox{ }  & & e^{-\delta \tilde{r}} \sum_{m=k_0}^{n} |a_{m}| |\gamma|^{m} \sum_{l=0}^{m-k_0} \delta^{-(l+1)} \frac{(m-k_0)!}{(m-k_0-l)!}  \tilde{r}^{m-l-k_0} \\
 &  & \phantom{mmm} \leq \delta^{-1} \sum_{m=k_0}^{n} |a_{m}| {r}^{m-k_0} |\gamma|^{k_0}  e^{-\delta \frac{{r}}{|\gamma|}}\sum_{l=0}^{m-k_0} \frac{(m-k_0)!}{(m-k_0-l)!}  \left( \delta \frac{{r}}{|\gamma|}\right)^{-l} \\
 & & \phantom{mmm}\leq |\gamma|^{n}\delta^{-(n-k_0 +1)} (n-k_0+1)! \\
 & & \phantom{mmmmmm} \sum_{m=k_0}^{n} |a_{m}| {r}^{m-n}   \sum_{l=0}^{m-k_0} \frac{(m-k_0)! (n-k_0-l)!}{(n-k_0+1)!(m-k_0-l)!} \frac{ \left( \delta \frac{{r}}{|\gamma|}\right)^{n-k_0-l} e^{-\delta \frac{{r}}{|\gamma|}}}{ (n-k_0-l)!}\\
 & & \phantom{mmm}\leq |\gamma|^{n+1}\delta^{-(n-k_0 +2)} (n-k_0+1)! \\
 & & \phantom{mmmmmm} \sum_{m=k_0}^{n-1} |a_{m}| {r}^{m-n}   \sum_{l=0}^{m-k_0} \frac{(m-k_0)! (n-k_0-l)!}{(n-k_0+1)!(m-k_0-l)!} \frac{1}{ \sqrt{2\pi(n-k_0-l)}}\\
  & & \phantom{mmmmmm} +  |a_{n}|   \left( \sum_{l=0}^{n-k_0-1} \frac{(n-k_0)! (n-k_0-l)!}{(n-k_0+1)!(n-k_0-l)!} \frac{1}{ \sqrt{2\pi(n-k_0-l)}} + \frac{1}{n-k_0+1}\right).\\
\end{eqnarray*}
Now observe that
$$ \frac{(m-k_0)! (n-k_0-l)!}{(n-k_0+1)!(m-k_0-l)!}  \leq \frac{1}{n-k_0+1},$$
and for $m<n$
$$ \sum_{l=0}^{m-k_0} \frac{1}{\sqrt{n-k_0 -l}} \leq \int_{n-m-1}^{n-k_0} \frac{1}{\sqrt{t}} dt =2 (\sqrt{n-k_0} - \sqrt{n-m-1})$$
and well as for $m=n$
$$ \sum_{l=0}^{n-k_0-1} \frac{1}{\sqrt{n-k_0 -l}} \leq \int_{0}^{n-k_0} \frac{1}{\sqrt{t}} dt = 2 \sqrt{n-k_0} .$$
Now observe that
\begin{eqnarray*} \frac{\sqrt{2}}{\sqrt{\pi}(n-k_0+1)}\sum_{m=k_0}^{n-1}  \sqrt{n-m-1} & \geq &
\frac{\sqrt{2}}{\sqrt{\pi}(n-k_0+1)}\int_{0}^{n-k_0} \sqrt{t} dt\\
&  =&  \frac{4\sqrt{2}}{3\sqrt{\pi}} \frac{(n-k_0)^{3/2}}{n+k_0 +1} \\
& \geq & \frac{1}{n-k_0+1}
\end{eqnarray*}
when ever
$$ n \geq \left(\frac{3 \sqrt{\pi}}{4\sqrt{2}}\right)^{2/3} + \  k_0$$
which is the case if $n \geq k_0 +2$. In the special case $n=k_0$ we see that
$$ \frac{1}{n-k_0+1} = 1 \leq \frac{2}{\sqrt{\pi(2n-2k_0+1)}}$$
and in the special case of $n= k_0 +1$ we get that
$$ \sum_{l=0}^{n-k_0} \frac{(n-k_0)! (n-k_0-l)!}{(n-k_0+1)!(n-k_0-l)!} \frac{ \left( \delta \frac{{r}}{|\gamma|}\right)^{n-k_0-l} e^{-\delta \frac{{r}}{|\gamma|}}}{ (n-k_0-l)!} = \frac{1}{2} (e^{-1} + e^{-\delta \frac{{r}}{|\gamma|}} )\leq \frac{1}{2} (e^{-1} +1).$$
But then we can continue to obtain for $n \neq k_0 +1$
 \begin{eqnarray*} 
\mbox{ }  & & e^{-\delta \tilde{r}} \sum_{m=k_0}^{n} |a_{m}| |\gamma|^{m} \sum_{l=0}^{m-k_0} \delta^{-(l+1)} \frac{(m-k_0)!}{(m-k_0-l)!}  \tilde{r}^{m-l-k_0} \\
& & \phantom{mmm}\leq |\gamma|^{n}\delta^{-(n-k_0 +1)} (n-k_0+1)! \\
 & & \phantom{mmmmmm} \frac{2{r}^{-n}}{\sqrt{\pi (2n-2k_0+1)}} \sum_{m=k_0}^{n} |a_{m}| {r}^{m} \\
&  & \phantom{mmm} \leq \frac{2}{\sqrt{\pi (2n-2k_0+1)}}{r}^{-n}\kappa_{k_0,n}(r) |\gamma|^{n}\delta^{-(n-k_0 +1)} (n-k_0+1)! 
\end{eqnarray*}
and in the case $n=k_0+1$
\begin{eqnarray*} 
\mbox{ }  & & e^{-\delta \tilde{r}} \sum_{m=k_0}^{n} |a_{m}| |\gamma|^{m} \sum_{l=0}^{m-k_0} \delta^{-(l+1)} \frac{(m-k_0)!}{(m-k_0-l)!}  \tilde{r}^{m-l-k_0} \\
&  & \phantom{mmm} \leq \frac{1}{2} (e^{-1} +1) {r}^{-n}\kappa_{k_0,n}(r) |\gamma|^{n}\delta^{-(n-k_0 +1)} (n-k_0+1)! 
\end{eqnarray*}
By appealing to Lemma \ref{0phiGk_0}, again using $\delta_1 > \delta$ and collecting terms, we get for $n=k_0-1$ that
$$
\int_{\tilde{\Gamma}_w} |K(w,z)| |\varphi_n(\gamma z)| |dz| \leq  \frac{dc_w}{b} |\gamma|^{n+1}   \delta^{-1} \inf_{r\in (0,R_f)} \frac{\left( \tilde{c}_\gamma {r}^{-n_0} + \kappa({r})\right) }{r^{n+1}},
$$
for $n=k_0$
\begin{eqnarray*}
\int_{\tilde{\Gamma}_w} |K(w,z)| |\varphi_n(\gamma z)| |dz| \leq  \frac{dc_w}{b\sqrt{2\pi}} |\gamma|^{n}   \delta^{-1} \phantom{yyyy} \\
\inf_{r\in (0,R_f)} \frac{\left( \tilde{c}_\gamma {r}^{-n_0} + \sqrt{2\pi} \kappa_{-n_0,k_0}({r}) + 2 \sqrt{2} \kappa_{k_0,k_0}(r) + \kappa_{k_0+1}(r)\right) }{r^{n}},
\end{eqnarray*}
for $n=k_0+1$
\begin{eqnarray*}
\int_{\tilde{\Gamma}_w} |K(w,z)| |\varphi_n(\gamma z)| |dz| \leq  \frac{dc_w}{b\sqrt{2\pi(n-k_0+1)}} |\gamma|^{n} \delta^{-(n-k_0+1)}  (n-k_0+1)! \phantom{yyyy} \\
\inf_{r\in (0,R_f)} \frac{\left( \tilde{c}_\gamma {r}^{-n_0}+ \frac{\kappa_{-n_0,k_0-1}({r})}{n-k_0+1}  +  \sqrt{\pi}(e^{-1} +1)\kappa_{k_0,n}(r) +  \kappa_{n}({r})\right)}{r^{n}}
\end{eqnarray*}
and for $n>k_0+1$
\begin{eqnarray*}
\int_{\tilde{\Gamma}_w} |K(w,z)| |\varphi_n(\gamma z)| |dz| \leq  \frac{dc_w}{b\sqrt{2\pi(n-k_0+1)}} |\gamma|^{n} \delta^{-(n-k_0+1)}  (n-k_0+1)! \phantom{yyyy} \\
\inf_{r\in (0,R_f)} \frac{\left( \tilde{c}_\gamma {r}^{-n_0}+ \frac{\kappa_{-n_0,k_0-1}({r})}{n-k_0+1}  +  2\sqrt{2}\frac{\sqrt{n-k_0+1}}{\sqrt{2n-2k_0+1}}\kappa_{k_0,n}(r) +  \kappa_{n}({r})\right)}{r^{n}}.
\end{eqnarray*}

 It might however be hard to use these estimates for $n\geq k_0$ in practice, since it requires control of $ \kappa_{-n_0, k_0-1}(r)$, $\kappa_{k_0, n}(r)$ and $\kappa_{n+1}(r)$ separately and the infimum involves the $\gamma$-dependent $\tilde{c}_\gamma$. So to obtain an estimate where the constant only involves $\kappa(r)$ and is completely $\gamma$ independent we have the slightly weaker estimate 
\begin{eqnarray*}
\int_{\tilde{\Gamma}_w} |K(w,z)| |\varphi_n(\gamma z)| |dz| \leq  \frac{2dc_w\tilde{c}_\gamma}{b\sqrt{\pi(1+ 2(n-k_0))}} |\gamma|^{n} \delta^{-(n-k_0+1)}  (n-k_0+1)! \phantom{yyyy} \\
\inf_{r\in (0,R_f)} \frac{\left( c {r}^{-n_0}+  \kappa({r})\right)}{r^{n}}
\end{eqnarray*}
which concludes the proof of Theorem \ref{Thm1}.
\eproof

We can also establish estimates for $n<k_0-1$, but under different and further assumptions. To this end fix an $r\in (0, R_f)$.

First of all, we will assume that there exist a continuous deformation $\overline{\Gamma}_w$, $w\in W$ of the contour $\Gamma_w$ such that 
$$\overline{\Gamma}_w \subset \bbC - (P_{K_w} \cup P_{f_\gamma}), \ \ \ \forall (w,\gamma) \in W\times (U\cap D(0, {r}))$$
and along which (\ref{cos}), (\ref{Mainestimateinfty}) and (\ref{Mainestimateinftyf}) holds.

Further, set
$$ \rho_{w} = \inf \{|z| \mid z\in \overline{\Gamma}_w\}$$
and assume it is positive 
\begin{equation}\label{rpos}
 \rho_w >0
 \end{equation}
for all $w \in W$. Introduce the constant
$$\overline{C}_n({r},\rho_w) =  \frac{ dc{r}^{-n_0} + d \max({r}^{-n_0}, \rho_w^{n-k_0+1})\kappa({r})}{b r^{n+1}}.$$

\begin{theorem}\label{Thm1.5}
For $n_0\leq n< k_0-1$ and any positive ${r}<R_f$  we have that
\begin{equation}\label{MENs}
\mid A^\gamma_{K,\Gamma}(f)(w)  -  \tilde A^{\gamma,n}_{K,\Gamma}(f)(w) \mid \leq \overline{C}_n({r},\rho_w) c_w \tilde{c}_\gamma|\gamma|^{n+1}   \delta^{-1} 
\end{equation}
for all $(w,\gamma) \in W\times (U\cap D(0, {r}))$.
\end{theorem}

In order to prove this theorem we first extend the above Lemmas. We make the same splitting of the contour as before
$$ \overline{\Gamma}^0_w( \tilde{r}) = \overline{\Gamma}_w \cap D(0,\tilde{r}), \ \ \  \overline{\Gamma}^{\infty}_w( \tilde{r}) = \overline{\Gamma}_w \cap (\bbC - D(0,\tilde{r})).$$

We now observe that the statement and proof of the first lemma is exactly the same, except we have to replace $\tilde{\Gamma}^0_w(\tilde{r})$ with $\overline{\Gamma}^0_w(\tilde{r})$ . The proof of the second lemma is also exactly the same and thus we immediately obtain.

\begin{lemma}\label{inftyznbar}
If $n < k_0-1$ then
$$
\sum_{m=-n_0}^{n} |a_m| \int_{\overline{\Gamma}_w^{\infty}(\tilde{r})} |K(w,z)| |\gamma z|^m |dz|  \leq  \frac{d c_w }{b}\delta_1^{-1} e^{-\delta_1 \tilde{r}}\sum_{m=-n_0}^{n} |a_{m}| |\gamma|^{m}  \tilde{r}^{m-k_0}$$
\end{lemma}

For the third lemma, we have in the case $n<k_0-1$.
\begin{lemma}\label{0phiSbar}
Provided that ${r} = \tilde{r}|\gamma|$ satisfies ${r}<R_f$, we have for $n <k_0-1$
\begin{equation*}%\label{Near0EN} 
\int_{\overline{\Gamma}_w^{0}(\tilde{r})} |K(w,z)| |\varphi_n(\gamma z)| |dz| \leq   \frac{d c_w  \kappa_{n+1}({r}){r}^{-n-1} \rho_w^{n-k_0+1}}{b}|\gamma|^{n+1} \delta_1^{-1}.
\end{equation*}
\end{lemma}

\proof
From the first part of the proof of Lemma \ref{0phiGk_0} we conclude that
$$
\int_{\overline{\Gamma}_w^{0}(\tilde{r})} |K(w,z)| |\varphi_n(\gamma z)| |dz|  \leq  \frac{d c_w  \kappa_{n+1}({r}){r}^{-n-1}}{b} |\gamma|^{n+1}  \int_{\rho_m}^{\tilde{r}} e^{-\delta_1 s} s^{n-k_0+1} |ds|.
$$
Thus for $n <k_0-1$ we get that
$$\int_{\overline{\Gamma}_w^{0}(\tilde{r})} |K(w,z)| |\varphi_n(\gamma z)| |dz| \leq \frac{d  c_w }{b} \kappa_{n+1}({r}){r}^{-n-1} \rho_w^{n-k_0+1} |\gamma|^{n+1} \delta_1^{-1} $$
\eproof

We now move on to proving Theorem \ref{Thm1.5}.

\proof
We assume that $|\gamma| \leq {r}$ and $n < k_0-1$.

With reference to Lemma \ref{inftyf}, we observe that 
 \begin{eqnarray*} 
\delta^{-1} |\gamma|^{-n_0} e^{-\delta \tilde{r}} \tilde{r}^{-(k_0+n_0)}  & \leq & |\gamma|^{n+1} (|\gamma|^{k_0-n-1} r^{-(k_0-n-1)}) {r}^{-n-n_0-1} \delta^{-1}\\
&\leq &  |\gamma|^{n+1}   {r}^{-n_0} r^{n-1} \delta^{-1}.
\end{eqnarray*}

Referring to Lemma \ref{inftyznbar}, using $\delta_1 > \delta$,
\begin{eqnarray*}
\delta^{-1}e^{-\delta \tilde{r}}\sum_{m=-n_0}^{n} |a_{m}| |\gamma|^{m}  \left(\frac{{r}}{|\gamma|}\right)^{m-k_0} & \leq &\delta^{-1} |\gamma|^{n+1} \kappa_{-n_0,n}({r})  (|\gamma|^{-n-1}{r}^{n+1}) r^{-k_0-n-1}\\
& \leq & \delta^{-1} |\gamma|^{n+1} \kappa_{-n_0,n}({r})   r^{-k_0-n-1} .
\end{eqnarray*}
Refering to Lemma \ref{0phiSbar}, using $\delta_1 >\delta$ and collecting terms we get that
\begin{equation*}%\label{MENs}
\mid A^\gamma_{K,\Gamma}(f)(w)  -  \tilde A^{\gamma,n}_{K,\Gamma}(f)(w) \mid \leq  \frac{d c_w |\gamma|^{n+1} }{b\delta r^{n+1}}   \left( \tilde{c}_\gamma r^{-n_0} + r^{-k_0} \kappa_{-n_0,n}(r) + \rho_w^{n-k_0+1} \kappa_{n+1}(r)\right)
\end{equation*}
from which the theorem follows.
\eproof

We can relax the condition on $\overline{\Gamma}_w$, by first of allowing it to depend on $\gamma$,  but we must of course have  
$$ \overline{\Gamma}_{w,\gamma} \subset \bbC - (P_{K_w} - P_{f_\gamma}), \ \ \ \forall (w,\gamma) \in W \times U.$$
Second, we restricting conditions ({\ref{cos}) to the part of $\overline{\Gamma}_w$ outside $D(0,\tilde{r})$, e.g. we assume
\begin{equation}\label{cosinfty}
\forall (w,\gamma) \in W\times U \ \ \exists r_{w,\gamma} \in (0,R_f) \ :    \ \ |\cos(\Theta(z)) | \geq b \ \ \ \ \forall z \in \overline{\Gamma}_{w,\gamma}^\infty(\frac{r_{w,\gamma}}{|\gamma|}).
\end{equation}
We now let
$$
\rho_{w,\gamma} = \inf \{|z| \mid z\in \overline{\Gamma}_{w,\gamma}\}
$$
and assume $\rho_{w,\gamma} >0$. Further let $L(r_{w,\gamma}) $ be the length of $\overline{\Gamma}_{w,\gamma}^0(\frac{r_{w,\gamma}}{|\gamma|})$.

We then need a replacement for Lemma \ref{0phiGk_0}.
 
 \begin{lemma}\label{LL}
Under condition (\ref{cosinfty}), (\ref{Mainestimateinfty}) and (\ref{Mainestimateinftyf}) we have that if $n< k_0, $  then
 $$
 \int_{\overline{\Gamma}_{w,\gamma}^0(\tilde{r})} \vert K(w,z) \phi_n(\gamma z)\vert \vert dz \vert \leq c_w \kappa_{n+1}(r) L_\gamma(r) r^{-n-1} \rho_{w,\gamma}^{n-k_0+1} |\gamma|^{n+1} \delta_1^{-1}.
 $$
 and further if $n \geq k_0$ 
 $$
  \int_{\overline{\Gamma}_{w,\gamma}^0(\tilde{r})} \vert K(w,z) \phi_n(\gamma z)\vert \vert dz \vert \leq \frac{L(r_{w,\gamma})c_w\kappa_{n+1}(r) r^{-n-1}}{\sqrt{2\pi (n-k_0+1)}} |\gamma|^{n+1} \delta_1^{n-k_0+1} (n-k_0+1)!
  $$
 \end{lemma}

\proof
We compute assuming $n<k_0$
\begin{eqnarray*}
\int_{\overline{\Gamma}_{w,\gamma}^0(\tilde{r})} \vert K(w,z) \phi_n(\gamma z)\vert \vert dz \vert & \leq & \int_{\overline{\Gamma}_{w,\gamma}^0(\tilde{r})} \vert K(w,z)\vert \vert \tilde{\phi}_n(\gamma z)\vert \vert \gamma z\vert^{n+1} \vert dz \vert \\
& \leq & c_w\kappa_{n+1}(r_{w,\gamma}) r_{w,\gamma}^{-n-1} |\gamma|^{n+1} \int_{\overline{\Gamma}_{w,\gamma}^0(\tilde{r})}  e^{-\delta |z|} \vert z \vert^{n+1-k_0}\vert dz \vert \\
& \leq & c_w\kappa_{n+1}(r_{w,\gamma}) r_{w,\gamma}^{-n-1} |\gamma|^{n+1} e^{-\delta_1 r_{w,\gamma} } \rho_{w,\gamma}^{n+1-k_0} \int_{\overline{\Gamma}_{w,\gamma}^0(\tilde{r})}  \vert dz \vert \\
& \leq & c_w\kappa_{n+1}(r_{w,\gamma}) r_{w,\gamma}^{-n-1} |\gamma|^{n+1} e^{-\delta_1 r_{w,\gamma} } \rho_{w,\gamma}^{n+1-k_0} L_\gamma(r)\\
\end{eqnarray*}
If $n\geq k_0$ then
\begin{eqnarray*}
\int_{\overline{\Gamma}_{w,\gamma}^0(\tilde{r})} \vert K(w,z) \phi_n(\gamma z)\vert \vert dz \vert & \leq &  c_w\kappa_{n+1}(r_{w,\gamma}) r_{w,\gamma}^{-n-1} |\gamma|^{n+1} \int_{\overline{\Gamma}_{w, \gamma}^0(\tilde{r})}  e^{-\delta_1 |z|} \vert z \vert^{n+1-k_0}\vert dz \vert  \\
& \leq & \frac{L(r_{w,\gamma})c_w\kappa_{n+1}(r_{w,\gamma}) r_{w,\gamma}^{-n-1}}{\sqrt{2\pi (n-k_0+1)}} |\gamma|^{n+1} \delta_1^{n-k_0+1} (n-k_0+1)!
\end{eqnarray*}

\eproof

Following the proof of Theorem \ref{Thm1} line for line, and instead of appealing to Lemma \ref{0phiGk_0} one appeals to Lemma  \ref{LL}, we obtain the following, using again $\delta_1 >\delta$.

For $n<k_0$
$$\int_{\overline{\Gamma}_{w,\gamma}} |K(w,z)| |\varphi_n(\gamma z)| |dz|  \leq \frac{dc_w}{b} |\gamma|^n \left(\tilde{c}_\gamma r_{w,\gamma}^{-n_0} + \kappa_{-n_0,n-1}(r_{w,\gamma}) + L(r_{w,\gamma})\kappa_n(r_{w,\gamma})\right).$$
and for $n\geq k_0$ we have that
\begin{eqnarray*}
\int_{\overline{\Gamma}_{w,\gamma}} |K(w,z)| |\varphi_n(\gamma z)| |dz| \leq |\gamma|^{n} \delta^{n-k_0+1} (n-k_0+1)! \phantom{jjjjjjjjjjjjhhhhhhhhhhhjjjjj}\\
\phantom{jjjj}   \frac{d c_w\tilde{c}_\gamma \left( c r_{w,\gamma}^{-n_0} + \frac{\kappa_{-n_0,k_0-1}(r_{w,\gamma})}{n-k_0+1} + 2\sqrt{2} \frac{\sqrt{n-k_0+1}}{\sqrt{2n-2k_0+1}}\kappa_{k_0,n}(r_{w,\gamma}) + L(r_{w,\gamma}) \kappa_n(r_{w,\gamma}) \right)}{b \sqrt{2\pi (n-k_0+1)} r_{w,\gamma}^{n+1}}.
\end{eqnarray*}
Let for $n< k_0$
$$\overline{C}^{(n)}_{w,\gamma} =   r_{w,\gamma}^{-n-1} \left( \frac{d c}{b} r_{w,\gamma}^{-n_0} + \max\left(1,L(r_{w,\gamma})\right) \kappa(r_{w,\gamma})\right) $$%c''_n(r_{w,\gamma}),$$
and for $n\geq k_0$ 
$$\overline{C}^{(n)}_{w,\gamma} =  \frac{1}{\sqrt{\pi (n-k_0+1)}} r_{w,\gamma}^{-n} \left( \frac{d c}{b} r_{w,\gamma}^{-n_0} + \max\left(2\sqrt{2}, \frac{L(r_{w,\gamma})}{r_{w,\gamma}}\right) \kappa(r_{w,\gamma})\right). $$ %c''_n(r_{w,\gamma}),$$
%where 
%$$ c''_n(r)  = .$$
We have then proved the following theorem.
\begin{theorem}\label{Thm1.75}
Under condition (\ref{cosinfty}), (\ref{Mainestimateinfty}) and (\ref{Mainestimateinftyf})  we have an asymptotic expansion in the Poincare sense to all orders
\begin{equation}\label{gamma_expansion2}
A^\gamma_{K,\Gamma}(f) \sim \tilde A^\gamma_{K,\Gamma}(f).
\end{equation}
In fact we have the following estimates for all $n \geq k_0$ 
\begin{equation}\label{MEPsss}
\mid A^\gamma_{K,\Gamma}(f)(w)  -  \tilde A^{\gamma,n}_{K,\Gamma}(f)(w) \mid \leq \overline{C}^{(n)}_{w,\gamma} c_w\tilde{c}_\gamma |\gamma|^{n} \delta^{n-k_0+1} (n-k_0+1)!
\end{equation}
and further for $n<k_0$, we have that 
\begin{equation}\label{MEPssss}
\mid A^\gamma_{K,\Gamma}(f)(w)  -  \tilde A^{\gamma,n}_{K,\Gamma}(f)(w) \mid \leq \overline{C}^{(n)}_{w,\gamma} c_w\tilde{c}_\gamma |\gamma|^{n+1} \delta^{-1} 
\end{equation}
for all $(w,\gamma) \in W\times U$.
\end{theorem}

We observe that $\overline{C}^{(n)}_{w,\gamma}$ is independent of $\rho_{w,\gamma}$. Further if we have that there exist $b>0$ such that
$$ \vert \cos(\Theta(z)) \vert \geq b \ \ \ \forall z \in \overline{\Gamma}_{w,\gamma}^\infty\left(\frac{r}{|\gamma|}\right), \ \ \ \forall r \in (0,R_f), \ \ \ (w,\gamma) \in W\times U,$$
then we get that (\ref{MEPsss}) also holds with the following definition
$$ \overline{C}^{(n)}_{w,\gamma}  = \inf_{r\in(0,R_f)}  r^{-n} \left( \frac{d c}{b} r^{-n_0} + \max\left(2\sqrt{2}, \frac{L(r)}{r}\right) \kappa(r)\right).
$$
We remark that in this case $\overline{C}^{(n)}_{w,\gamma} $ is actually independent of $(w,\gamma) \in W\times U$.

\section{The simpler case of a single pole at zero.}
\label{Mild}

In this section we will establish an analog of Theorem \ref{Thm1} in the case where 
$$R_f = \infty.$$ 
We will still assume the setup from the introduction with the assumptions (\ref{cos}), (\ref{Mainestimateinfty}), however we will not need (\ref{Mainestimateinftyf}). But we will now make the following ekstra assumption. Recall that
$$ \tilde{\phi}_n(z) = \frac{\phi_n(z)}{z^{n+1}} = \sum_{m=n+1}^\infty a_m z^{m-n-1}.$$
For $n\geq k_0-1$, we will assume there exist a constant $\tilde{\delta}$ and
$$ C^{(n)} : W \rightarrow \bbR_+$$
such that
\begin{equation} \label{tpest}
  |\tilde{\phi}_n(\gamma z)| \leq C^{(n)}_w e^{\tilde{\delta} |z|}
\end{equation}
for all $z\in \tilde{\Gamma}_w$ and $(w,\gamma) \in W\times U$. Let now
$$ \delta = \delta_1 - \tilde{\delta}.$$
\begin{theorem}\label{Thm77}
Under the assumptions (\ref{cos}), (\ref{Mainestimateinfty}) and (\ref{tpest}), we have an asymptotic expansion in the Poincare sense to all orders
\begin{equation}\label{gamma_expansion2}
A^\gamma_{K,\Gamma}(f) \sim \tilde A^\gamma_{K,\Gamma}(f).
\end{equation}
In fact we have the following estimates for all $n \geq k_0-1$ 
\begin{equation}\label{MEPs}
\mid A^\gamma_{K,\Gamma}(f)(w)  -  \tilde A^{\gamma,n}_{K,\Gamma}(f)(w) \mid \leq C^{(n)}_w \frac{d c_w\tilde{c}_\gamma}{b} |\gamma|^{n+1} \delta^{n-k_0+2} (n-k_0+1)!
\end{equation}
for all $(w,\gamma) \in W\times U$.
\end{theorem}
\proof
The proof starts the same ways as the proof of Theorem \ref{Thm1}, so assume that $n\geq k_0-1$, we get that 
$$z\mapsto K(w,z) \varphi_n(\gamma z) \text{ is regular } \forall z\in D(0,R_K(w))$$
for all $w\in W$, where we recall that $R_K(w)$ is the minimal distance from zero to $P_{K_w} - \{0\}$. 
Thus we see that
$$ A^\gamma_{K,\Gamma}(f)(w)  -  \tilde A^{\gamma,n}_{K,\Gamma}(f)(w) = \int_{\tilde{\Gamma}_w} K(w,z) \varphi_n(\gamma z)dz$$
and therefore
\begin{eqnarray*}
  \mid A^\gamma_{K,\Gamma}(f)  -  \tilde A^{\gamma,n}_{K,\Gamma}(f) \mid &\leq& \int_{\tilde{\Gamma}_w} |K(w,z)| |\varphi_n(\gamma z)| |dz|\\
  &\leq &c_w \int_{\tilde{\Gamma}_w} e^{- \delta_1 |z|} |z|^{-k_0}|\tilde{\varphi}_n(\gamma z)| |\gamma z|^{n+1}|dz|\\
  &\leq &c_w C_w^{(n)} |\gamma|^{n+1}\int_{\tilde{\Gamma}_w} e^{- \delta|z|} |z|^{-k_0 + n+1}|dz|\\
  &\leq &c_w C_w^{(n)} \frac{d}{b} |\gamma|^{n+1}\int_{0}^\infty  e^{- \delta s} s^{n -k_0 +1}ds\\
  &=&c_w C_w^{(n)} \frac{d}{b} |\gamma|^{n+1} \delta^{n-k_0+2} (n-k_0+1)!
 \end{eqnarray*}

\eproof

We observe that we can under an extra assumption get a formula for the Laurent expansion in $\gamma$ at zero.

\begin{proposition}\label{PropConv}
Assume that
$$\limsup_{m\rightarrow \infty} (|a_m| (m-k_0)!)^{1/m} = R^{-1} < \infty$$
then the Laurent series for $g_\gamma(w)$ is 
$$ g_\gamma(w) = \sum_{m=-n_0}^\infty a_m h_m(w) \gamma^m,$$
which converges for $0<|\gamma| < \delta_1 R$.
\end{proposition}

\proof
We compute for $m\geq k_0$ that
\begin{eqnarray*}
  \mid h_m(w)\mid &\leq& \int_{\tilde{\Gamma}_w} |K(w,z)| |z|^m |dz|\\
  &\leq &c_w \int_{\tilde{\Gamma}_w} e^{- \delta_1 |z|} |z|^{m- k_0}|dz|\\
  &=&c_w  \frac{d}{b} \delta_1^{m-k_0+1} (m-k_0)!
 \end{eqnarray*}
The proposition follows immediately from this.

\eproof

We would like to also consider the case where we restrict the assumption ({\ref{cos}) as we have done in the previous section, e.g. we assume there exist 
$$ r : W \rightarrow \bbR_+$$
and a continuous deformation of $\Gamma_w$ to a contour 
$$\overline{\Gamma}_w \subset \bbC -(P_{K_w} \cup P_(f_\gamma))$$
for all $(w,\gamma) \in W\times U$, such that (\ref{cos}), (\ref{Mainestimateinfty}) and (\ref{tpest}) are satisfied along
$$ \overline{\Gamma}^\infty_w(r(w)) := \overline{\Gamma}_w \cap (\bbC - D(0,r(w))).$$
As before we also let $\overline{\Gamma}^0_w(r(w)) :=  \overline{\Gamma}_w \cap D(0,r(w))$ and we can now with these assumptions prove the following theorem.

\begin{theorem}
Under the assumptions of (\ref{cos}), (\ref{Mainestimateinfty}) and (\ref{tpest}) along $\overline{\Gamma}^\infty_w(r(w))$, we have an asymptotic expansion in the Poincare sense to all orders
\begin{equation}\label{gamma_expansion2}
A^\gamma_{K,\Gamma}(f) \sim \tilde A^\gamma_{K,\Gamma}(f).
\end{equation}
In fact there exist constants $\tilde{C}^{(n)}_w$ for all $n \geq - n_0$ and $w\in W$ such that
\begin{equation}\label{MEPs}
\mid A^\gamma_{K,\Gamma}(f)(w)  -  \tilde A^{\gamma,n}_{K,\Gamma}(f)(w) \mid \leq \tilde{C}^{(n)}_w  |\gamma|^{n+1} 
\end{equation}
for all $(w,\gamma) \in W\times U$.
\end{theorem}

\proof
We have that
$$ A^\gamma_{K,\Gamma}(f)(w)  -  \tilde A^{\gamma,n}_{K,\Gamma}(f)(w) = \int_{\tilde{\Gamma}_w} K(w,z) \varphi_n(\gamma z)dz$$
and therefore
\begin{eqnarray*}
  \mid A^\gamma_{K,\Gamma}(f)  -  \tilde A^{\gamma,n}_{K,\Gamma}(f) \mid &\leq& \int_{\overline{\Gamma}^0_w(r(w))} |K(w,z)| |\tilde{\varphi}_n(\gamma z)| |\gamma z|^{n+1} |dz|\\
  && \phantom{HH} + \int_{\overline{\Gamma}^\infty_w(r(w))} |K(w,z)||\tilde{\varphi}_n(\gamma z)| |\gamma z|^{n+1} |dz|\\
  &\leq &M_w L(r(w)) |\gamma|^{n+1}\\
  & & \phantom{HH} +c_w \int_{\overline{\Gamma}^\infty_w(r(w))} e^{- \delta_1 |z|} |z|^{-k_0}|\tilde{\varphi}_n(\gamma z)| |\gamma z|^{n+1}|dz|\\
  &\leq &M_w L(r(w)) |\gamma|^{n+1}\\
  & & \phantom{HH} +c_w C_w^{(n)} |\gamma|^{n+1}\int_{\overline{\Gamma}_w} e^{- \delta |z|} |z|^{-k_0 + n+1}|dz|\\
  &\leq &C_w^{(n)} M_w L(r(w)) |\gamma|^{n+1}\\
  & & \phantom{HH} +c_w C_w^{(n)} \frac{d}{b} |\gamma|^{n+1}\int_{0}^\infty  e^{- \delta s} s^{n -k_0 +1}ds\\
  &\leq &C_w^{(n)}  \left( M_w L(r(w))  + c_w \frac{d}{b} \delta^{n-k_0+2} (n-k_0+1)!\right) |\gamma|^{n+1} ,
 \end{eqnarray*}
where 
$$ M_w = \max_{z\in \overline{\Gamma}^0_w(r(w))} |K(w,z)| e^{\tilde{\delta} |z|} |z|^{n+1}$$
and $L(r(w))$ is the length of $\overline{\Gamma}^0_w(r(w))$.
\eproof

\section{Borel resummation}

\label{Borel}

We will now prove Theorem \ref{Thm2}.

\proof

We start by verifying that the integral expression for $B_w(\xi)$ is absolutely integrable
\begin{eqnarray*}
 |B_w(\xi)| &\leq & \int_{\overline{\Gamma}_{w}} |K(w,z)| \vert \calL_\theta^{-1}(\phi_z)(\xi)\vert |dz| \\
  &\leq & \int_{\overline{\Gamma}^0_{w}(r(w))} |K(w,z)| \vert \calL_\theta^{-1}(\phi_z)(\xi)\vert |dz| \\
 & & + \int_{\overline{\Gamma}^\infty_{w}(r(w))} e^{-\delta_1 |z|} |z|^{-k_0} \vert \calL_\theta^{-1}(\phi_z)(\xi)\vert |dz| \\
& \leq & \left( M_w L(r(w)) + \frac{d e^{-\delta_1 r(w)} r(w)^{-k_0}}{b}\right) \sum_{i} C_i e^{\alpha_\theta(w)|\xi|} |\xi |^{m_i}
\end{eqnarray*}
where
$$ M_w = \max_{z\in \overline{\Gamma}^0_{w}(r(w))} |K(w,z)|$$
and $L(r(w))$ is the length of $\overline{\Gamma}^0_{w}(r(w))$.
Letting
$$
\tilde{C}_i = \left( M_w L(r(w)) + \frac{d e^{-\delta_1 r(w)} r(w)^{-k_0}}{b}\right) C_i 
$$
we further compute
\begin{eqnarray*}
\int_{\bbR_+^\theta} \int_{\overline{\Gamma}_{w}} |e^{-\xi/\gamma} K(w,z) \calL_\theta^{-1}(\phi_z)(\xi)| |dz| |d\xi| &\leq &\sum_i \tilde{C}_i \int_{\bbR_+^\theta} |e^{-\xi/\gamma}| e^{\alpha_\theta(w)|\xi|} |\xi|^{m_i} |d\xi |\\
&\leq &  \sum_i \tilde{C}_i  \int_{\bbR_+^\theta} e^{-(\frac{1}{|\gamma|} \cos(\theta - \theta_\gamma) - \alpha_\theta(w))|\xi|} |\xi|^{m_i} |d\xi |\\
&< &  \infty
 \end{eqnarray*}
 since $\gamma \in U_{\theta,\alpha(w)}$
for all $ w \in W$.
Thus by Fubini's theorem we get that
\begin{eqnarray*}
\calL_\theta(B_w)(\gamma) & = &  \int_{\bbR_+^\theta} e^{-\gamma/\xi} \left(\int_{\tilde{\Gamma}_{w}} K(w,z)  \calL_\theta^{-1}(\phi_z)(\xi) dz\right) d\xi\\
& = &  \int_{\tilde{\Gamma}_{w}} K(w,z)   \left(\int_{\bbR_+^\theta} e^{-\xi/\gamma}\calL_\theta^{-1}(\phi_z)(\xi) d\xi\right) dz\\
 & = &\int_{\tilde{\Gamma}_{w}} K(w,z)  \phi_z(\gamma)  dz\\
 & = &g_\gamma(w)  -g^-_\gamma(w) .
\end{eqnarray*}
which establishes (\ref{Laplace-Borel}).

Let the Taylor series in zero of $B_w(\xi)$ be
$$ B_w(\xi) = \sum_{l=0}^\infty b_l(w) \xi^l,$$
which is uniform absolute convergent for all $\xi \in D(0,\tilde{r}_w)$, where $\tilde{r}_w$ is less than the radius of convergence. We recall that this implies that there exist a constant $C_L(w)$ such that
$$ \vert B_w(\xi) - \sum_{l=0}^{L-1} b_l(w) \xi^l \vert \leq C_L(w) |\xi|^{L}$$
for all $\xi \in D(0,\tilde{r}_w)$ and all $w\in W$.  
We claim that there exist constants $\tilde{C}_L(w)$  for all $w\in W$ such that
\begin{equation}\label{Lapest}
 \bigg\vert \calL_\theta(B_w)(\gamma)  - \sum_{l=0}^{L-1} b_l(w) l! \gamma^{l+1} \bigg\vert \leq  \tilde{C}_L(w) |\gamma|^{L+1}
 \end{equation}
for all $\gamma \in U_{\theta,\alpha(w)}$.
We have of course that
\begin{eqnarray*}\calL_\theta(B_w)(\gamma) - \sum_{l=0}^{L-1} b_l(w) l! \gamma^{l+1} &=& \int_{e^{i\theta}[0,\tilde{r}_w]} e^{-\xi/\gamma} \left(B_w(\xi)  - \sum_{l=0}^{L-1} b_l(w) \xi^l \right)d\xi \\
& & + \int_{e^{i\theta}[\tilde{r}_w, \infty)} e^{-\xi/\gamma}  \left(B_w(\xi)  - \sum_{l=0}^{L-1} b_l(w) \xi^l \right) d\xi.
\end{eqnarray*}
The first of these integrals we estimate as follows
\begin{eqnarray}
\bigg\vert \int_{e^{i\theta}[0,\tilde{r}_w]} e^{-\xi/\gamma} \left(B_w(\xi)  - \sum_{l=0}^{L-1} b_l(w) \xi^l \right) d\xi \bigg\vert \nonumber\\
\leq   C_L(w) \int_{e^{i\theta}[0,\tilde{r}_w]} \vert e^{-\xi/\gamma} \vert  |\xi|^{L} \vert d\xi \vert  \nonumber\\
\leq   C_L(w) |\gamma|^{L+1}. \label{1}
\end{eqnarray}
For all $\gamma \in U_{\theta, \alpha(w)}$
\begin{eqnarray*}
\bigg\vert  \int_{e^{i\theta}[\tilde{r}_w, \infty)} e^{-\xi/\gamma}  B_w(\xi)  d\xi \bigg\vert \\
& \leq & \int_{\tilde{r}_w}^\infty e^{ - \frac{1}{|\gamma|}\cos(\theta-\theta_\gamma) t} |B_w(e^{i\theta}t)| dt\\
& \leq & \sum_i \tilde{C}_i\int_{\tilde{r}_w}^\infty e^{ - (\frac{1}{|\gamma|}\cos(\theta-\theta_\gamma) - \alpha_\theta(w)) t} |t|^{m_i} dt\\
& \leq &  e^ {-\tilde{\delta}(\gamma)\tilde{r}_w}\sum_i \tilde{C}_i \sum_{l=0}^{m_i} \tilde{\delta}(\gamma)^{-l-1} \tilde{r}_w^l \frac{m_i!}{(m_i-l)!}
\end{eqnarray*}
where
$$\tilde{\delta}(\gamma) = \frac{1}{|\gamma|} \cos(\theta - \theta_\gamma) - \alpha_\theta(w) >0.$$
The fact that $\gamma \in U_{\theta}$ means that 
$$ \cos(\theta - \theta_\gamma) >0.$$
We observe that
\begin{eqnarray*}
 e^{-\frac{1}{|\gamma|} \cos(\theta - \theta_\gamma)} &= &\frac{e^{-\frac{1}{|\gamma|} \cos(\theta - \theta_\gamma)} \cos^{L+1}(\theta - \theta_\gamma) |\gamma|^{L+1}}{|\gamma|^{L+1} \cos^{L+1}(\theta - \theta_\gamma)}\\
&\leq & \frac{(L+1)! }{\sqrt{2\pi (L+1)}} \frac{|\gamma|^{L+1}}{\cos^{L+1}(\theta - \theta_\gamma) }.
\end{eqnarray*}
We now make the extra assumption on $\gamma$, namely that there exist a positive $c_\theta(w)$ depending only on $w$ such that
$$
\bigg\vert \frac{\cos (\theta - \theta_\gamma)}{|\gamma|} - \alpha_\theta(w)\bigg\vert >c_\theta(w)
$$
which is clearly satisfies for $|\gamma |$ small enough. 
Then we conclude that
\begin{equation}\label{2}
\bigg\vert  \int_{e^{i\theta}[\tilde{r}_w, \infty)} e^{-\xi/\gamma}  B_w(\xi)  d\xi \bigg\vert  \leq \frac{  e^{\alpha_\theta(w) \tilde{r}_w} (L+1)! }{ \sqrt{2\pi (L+1)}}\sum_i \tilde{C}_i\left( \sum_{l=0}^{m_i} c_\theta(w)^{-l-1} \tilde{r}_w^l \frac{m_i!}{(m_i-l)!}\right)
\frac{|\gamma|^{L+1}}{\cos^{L+1}(\theta - \theta_\gamma)}.
\end{equation}
Now we consider
\begin{eqnarray}
\bigg\vert \int_{e^{i\theta}[\tilde{r}_w, \infty)} e^{-\xi/\gamma}  \sum_{l=0}^{L-1} b_l(w) \xi^l d\xi \bigg\vert & \leq & \sum_{l=0}^{L-1} |b_l(w)| \int_{\tilde{r}_w}^\infty e^{- \frac{1}{|\gamma|}\cos(\theta- \theta_\gamma) t} t^l dt\nonumber \\
& \leq &  \sum_{l=0}^{L-1} |b_l(w)|\sum_{l'=0}^l \frac{|\gamma|^{l'+1} l!}{\cos^{l'+1}(\theta - \theta_\gamma)} \frac{ e^{-\frac{1}{|\gamma|}\cos(\theta - \theta_\gamma) \tilde{r}_w}}{(l-l')!}\nonumber \\
& \leq &  \left( \sum_{l=0}^{L-1} |b_l(w)|\frac{l!}{\cos^{L+1}(\theta - \theta_\gamma)}\right. \nonumber\\
& & \phantom{HH}\left.\sum_{l'=0}^l \frac{ (L-l')!}{(l-l')!\sqrt{2\pi(L-l')}\tilde{r}_w^{L-l'} }\right) |\gamma|^{L+1}. \label{3}
\end{eqnarray}
Combining (\ref{1}), (\ref{2}) and (\ref{3}) we conclude (\ref{Lapest}). Recall now the global estimates of Theorem \ref{Thm1} in $U$  on the asymptotic expansion of $g_\gamma(w)  -g^-_\gamma(w)$ in $\gamma$, which shows that the asymptotic expansion of $g_\gamma(w)  -g^-_\gamma(w)$ is $g_\gamma^+(\xi)$.  But then by (\ref{Laplace-Borel}) we must have that
 $$ \sum_{l=0}^\infty b_l(w) \xi^l = \mathcal{B}(g^+_\gamma)(\xi)$$ 
 for all $\xi \in D(0,\tilde{r}_w)$ and all $w\in W$, establishing (\ref{Borel-resum}).

\eproof

\section{Resurgence and Stokes coefficients}

\label{Stokes}

The resurgence properties of $B_w$ is really nicely understood using our formula
$$ B_w(\xi) = \int_{\overline{\Gamma}_{w}} K(w,z) \calL_\theta^{-1}(\phi_z)(\xi) dz.$$
Further the Stokes lines and the Stokes jumps can also be derived from the these expressions by considering the pole structure of both $\phi_z$ and further of $B_w$ it self as is clear from Theorem \ref{Thm3} and \ref{Thm5} which we no establish starting with Theorem \ref{Thm3}.

\proof
We recall the setting of Theorem \ref{Thm3} and the assumptions of this theorem.

By assumption $\calL_\theta(B_w)$ is well defined for $\theta \in (\theta_j, \theta_j\pm \epsilon)$. In general we see that $\calL_\theta(B_w)$ will be constant on sectors where $B_w$ has no poles  and decays at infinity sufficiently fast so as to make $\calL_\theta(B_w)$  well defined. However, if we consider a line $\bbR^{\theta}_+$ which hits the poles of $B_w$, then $\calL_\theta(B_w)$ is not well-defined along these lines. These are the directions  
$\theta_J = \{ \theta_j \mid j \in J_w\}$ index by set $J_w$, $w\in W$ as introduced in the introduction. 

We now consider 
$$ \Gamma'_n = e^{i(\theta_j-\epsilon)} [0, R_n] \cup \Gamma_n \cup e^{i(\theta_j+\epsilon)} [R_n, 0].$$
By the residue theorem we have that
$$ \int_{\Gamma'_n} e^{-\xi/\gamma} B_w(\xi) d\xi = 2 \pi i \sum_{p\in P_{B_w} \cap \bbR^{\theta_j}_+\cap D(0,R_n)} \Res_{\xi=p} ( e^{-\xi/\gamma}B_w(\xi)).$$
The theorem now follows directly by taking the limit as $n\rightarrow \infty$ and using that
$$ 
 \lim_{n\rightarrow \infty}  \bigg\vert \int_{\Gamma'_n} e^{-\xi/\gamma}B_w(\xi)) d\xi \bigg\vert = 4 \pi  \lim_{ n\rightarrow \infty}  R_n \sup_{\xi \in \Gamma'_n} \vert e^{-\xi/\gamma}B_w(\xi)) \vert  =  0.
$$
\eproof

We further prove Theorem \ref{Thm5} as follows.

\proof Following the same line of arguments as in the proof of Theorem \ref{Thm2} and using the assumptions on $\psi_z$ of Theorem \ref{Thm5}, we see that the integral
$$ \Psi_{w}^{(2)}(\xi)  := \int_{\overline{\Gamma}_{w}} K(w,z) \calL_\theta^{-1}(\psi_z)(\xi) dz$$
is absolutely convergent and defines  $\Psi_{w}^{(2)} \in \calO(V_{w,\theta})$ for all $w\in W$. 
Now introduce
$$
 \tilde{\Phi}_z(\xi)  =  -\sum_{p\in P_{\phi}}  \sum_{m=1}^{n_\varphi} \frac{b_{p,m}}{p^m}  \sum_{l=0}^{m-1}  {{m}\choose{l}}(-1)^l \sum_{l'=0}^l  {{l}\choose{l'}} \left(\frac{z}{p}\right)^{m-l+l'}  \frac{e^{\frac{z\xi}{p}}\xi^{m-l+l'-1}}{(m-l+l'-1)!} 
$$
We have that
\begin{eqnarray}
| \tilde{\Phi}_z(\xi)|  & \leq & \sum_{p\in P_{\phi}}  \sum_{m=1}^{n_\varphi} \frac{|b_{p,m}|}{|p|^m}  \sum_{l=0}^{m-1}  {{m}\choose{l}}\sum_{l'=0}^l  {{l}\choose{l'}} \left|\frac{z}{p}\right|^{m-l+l'} \nonumber \\
 & & \phantom{HH} \frac{e^{\frac{|z| |\xi |}{|p|}\cos(\theta_z + \theta - \theta_p)}|\xi|^{m-l+l'-1}}{(m-l+l'-1)!} \nonumber\\
 & \leq &    \sum_{m=1}^{n_\varphi}  \sum_{l=0}^{m-1}  {{m}\choose{l}}\sum_{l'=0}^l  {{l}\choose{l'}} \left|z\right|^{m-l+l'} \nonumber \\
 & & \phantom{HH} \frac{e^{\alpha_\theta(w) |\xi|}|\xi|^{m-l+l'-1}}{(m-l+l'-1)!} \sum_{p\in P_{\phi}}  \frac{|b_{p,m}|}{|p|^{2m-l+l'} } \label{estPhi}
 \end{eqnarray}
Thus the series expression for $\tilde{\Phi}$ is absolutely convergent by (\ref{Asump1}). We continue to estimate
\begin{eqnarray*}
\bigg\vert \int_{\bbR_+^\theta} e^{-\xi/\gamma} \tilde{\Phi}_z(\xi) d\xi \bigg\vert & \leq &   \sum_{m=1}^{n_\varphi}  \sum_{l=0}^{m-1}  {{m}\choose{l}}\sum_{l'=0}^l  {{l}\choose{l'}} \left|z\right|^{m-l+l'}  \\
& & \left(\int_{\bbR_+^\theta} e^{ - \left(\frac{1}{|\gamma|}\cos(\theta-\theta_\gamma) - \alpha_\theta(w) \right) | \xi|}|\xi|^{m-l+l'-1} |d\xi|\right) \sum_{p\in P_{\phi}}  \frac{|b_{p,m}|}{|p|^{2m-l+l'} } \\
& \leq &   \sum_{m=1}^{n_\varphi}  \sum_{l=0}^{m-1}  {{m}\choose{l}}\sum_{l'=0}^l  {{l}\choose{l'}} \left|z\right|^{m-l+l'}  \\
& &\frac{1}{\left(\frac{1}{|\gamma|}\cos(\theta-\theta_\gamma) - \alpha_\theta(w) \right)^{m-l+l'}} \sum_{p\in P_{\phi}}  \frac{|b_{p,m}|}{|p|^{2m-l+l'} } \\
\end{eqnarray*}
using (\ref{cosrel}). Thus we can use Fubini's theorem to compute
\begin{eqnarray*}
\int_{\bbR_+^\theta} e^{-\xi/\gamma} \tilde{\Phi}_z(\xi) d\xi & = &  -\sum_{p\in P_{\phi}}  \sum_{m=1}^{n_\varphi} \frac{b_{p,m}}{p^m}  \sum_{l=0}^{m-1}  {{m}\choose{l}}(-1)^l \sum_{l'=0}^l  {{l}\choose{l'}} \left(\frac{z}{p}\right)^{m-l+l'}   \\
& & \int_{\bbR_+^\theta} e^{ - \xi/\gamma} e^{\frac{z\xi}{p}}\xi^{m-l+l'-1} d\xi \\
& = & \tilde{\phi}'_z(\gamma)
\end{eqnarray*}
for all $\gamma \in U_{\theta, \alpha(w)}\cap O_{\epsilon,z}.$ Hence we conclude that
$$ \calL_\theta^{-1}(\tilde{\phi}'_z) (\xi) = \tilde{\Phi}_z(\xi), \ \ \ \forall \xi \in V_{w,\theta}.$$
We observe by (\ref{estPhi}) that condition (\ref{intest}) in Theorem \ref{Thm2} is satisfied.
We obtain further the following estimate
\begin{eqnarray}\label{C1est}
\int_{\overline{\Gamma}_{w}} |K(w,z)| |e^{\frac{z\xi}{p}}| \left|\frac{z}{p}\right|^{m-l+l'}  |dz| & \leq & \frac{e^{\alpha_\theta(w) |\xi|}}{|p|^{m-l+l'}} \int_{\overline{\Gamma}_{w}} |K(w,z)| \left|z\right|^{m-l+l'}  |dz| \nonumber\\
& \leq &  \frac{e^{\alpha_\theta(w) |\xi|}}{|p|^{m-l+l'}} \int_{\overline{\Gamma}^0_{w}(r(w))} |K(w,z)|\left|z\right|^{m-l+l'}  |dz| \nonumber\\
& &  + \frac{e^{\alpha_\theta(w) |\xi|}\rho_w^{-k_0}}{|p|^{m-l+l'}} \int_{\overline{\Gamma}^\infty_{w}(r(w))} e^{-\delta_1 |z|} \left|z\right|^{m-l+l'}  |dz| \nonumber\\
& \leq &  \frac{e^{\alpha_\theta(w) |\xi|}\rho_w^{-k_0}}{|p|^{m-l+l'}} M_w L(r(w)) r(w)^{m-l+l'} \\
& &  + \frac{de^{\alpha_\theta(w) |\xi|}\rho_w^{-k_0}}{b|p|^{m-l+l'}} \sum_{m'= 0}^{m-l+l'} \delta_1^{-m'-1} \frac{(m-l+l')!}{m'!} r(w)^{m-l+l' - m'}, \nonumber
\end{eqnarray}
where
$$ M_w = \max_{z \in \overline{\Gamma}_w^0(r(w))} |K(w,z)|.$$
Thus 
\begin{eqnarray*}
\Psi_{w}^{(1)}(\xi) &:= & -\sum_{p\in P_{\phi}}  \sum_{m=1}^{n_\varphi} \frac{b_{p,m}}{p^m}  \nonumber \sum_{l=0}^{m-1}  {{m}\choose{l}}(-1)^l \sum_{l'=0}^l  {{l}\choose{l'}}  \\
& & \phantom{jjjjj} \left( \int_{\overline{\Gamma}_{w}} K(w,z) e^{\frac{z\xi}{p}} \left(\frac{z}{p}\right)^{m-l+l'}  dz \right) \frac{\xi^{m-l+l'-1}}{(m-l+l'-1)!} \nonumber
\end{eqnarray*}
is absolutely convergent by (\ref{Asump1}) and  gives a well defined holomorphic function $\Psi^{(1)}_w \in \calO(V_\theta)$ for all $z \in \overline{\Gamma}_w$ and all $w\in W$.
We now get immediately  from Fubini's theorem that
$$ \Psi_{w}^{(1)}(\xi)  = \int_{\overline{\Gamma}_w} K(w,z) \calL_\theta^{-1}(\tilde{\phi}'_z)(\xi) dz \ \ \forall \xi \in V_{\theta,w}.$$
But then 
$$
B_w(\xi) = \int_{\overline{\Gamma}_w} K(w,z) \calL_\theta^{-1}(\tilde{\phi}'_z + \psi_x)(\xi) dz = \Psi_{w}^{(1)}(\xi)  + \Psi_{w}^{(2)}(\xi).
$$
for all $\xi \in V_{\theta,w}$, which concludes the proof.

\eproof

\section{The examples of Faddeev's quantum dilogarithm $S_\gamma$}

\label{Faddeev}

We let 
$$ K_F(w,z) = \frac{e^{w z}}{\sinh(\pi z) z} \text{ and } f_F(z) = \frac{1}{\sinh(z)}$$
and 
$$ S_\gamma(w) = \exp(\frac14 g^F_\gamma(w)) \ \ \ \forall (w,\gamma) \in \widetilde{W}^F\times \widetilde{U}^F$$
where
$$ g_\gamma^F(w) = \int_{\bbR+i\epsilon} \frac{e^{wz}}{\sinh(\pi z) z} \frac{1}{\sinh(\gamma z)} dz$$
and
$$ \widetilde{W}^F = \{ w\in \bbC \mid |\Re(w)| < \pi  + |\Re(\gamma)|\}, \ \ \ \widetilde{U}^F = \{ \gamma \in \bbC \mid \Re(\gamma) >0\}$$
where $S_\gamma \in \calM(\bbC)$ is Faddeev's quantum dilogarithm.

We observe that $k_0 = 2$ and $n_0 =1$ and further that $R_{K_F} = 1$ and $R_{f_F} = \pi$.
We recall that the Laurent series for $f_F$ convergent in $D(0,\pi)$ is 
$$ f_F(z) = \sum_{m=0}^\infty \frac{2 (1 - 2^{2m-1}) B_{2m}}{(2m)!} z^{2m-1} \ \ \ \forall z \in D(0,\pi),$$
where $B_{2m}$ is the $2m$'th Bernoulli number.
In anticipation of the asymptotic expansion we let
$$ P_\gamma^{2n}(w) = \frac{1}{4} \sum_{m=0}^n \frac{2 (1 - 2^{2m-1}) B_{2m}}{(2m)!} \gamma^{2m-1}\int_{\bbR + i\epsilon} \frac{e^{wz}}{\sinh(\pi z)z}  z^{2m-1} dz.$$
But starting with (for a proof of this formula see e.g.  \cite{AH})
$$ \frac{1}{2i} \Li_2(-e^{iw}) = \frac{1}{4} \int_{\bbR + i \epsilon}  \frac{e^{wz}}{\sinh(\pi z) z^2}dz$$
and differentiating $2m$ times in $w$ we obtain the $m$'th coefficient
$$ \frac{1}{2i} \left(\frac{\partial}{\partial w}\right)^{2m} \Li_2(-e^{iw}) = \int_{\bbR + i\epsilon} \frac{e^{wz}}{\sinh(\pi z)z}  z^{2m-1} dz$$
and so
$$ P_\gamma^{2n}(w) = \frac{1}{2i} \sum_{m=0}^n \frac{2 (1 - 2^{2m-1}) B_{2m}}{(2m)!} \gamma^{2m-1}\left(\frac{\partial}{\partial w}\right)^{2m} \Li_2(-e^{iw}) .$$

Since we will take $\tilde{\Gamma}^{F}_w = \tilde{\Gamma}^F_{\tilde{\theta}}$ to be the real line rotated by the angle $\tilde{\theta}$, $d^F=2$ and for $\tilde{\theta} \in (-\frac{\pi}{2},\frac{\pi}{2})$ we let
$$\tilde{\Gamma}^F_{\tilde{\theta}} = e^{i\tilde{\theta}}(-\infty, 0]\cup e^{i\tilde{\theta}} [0,\infty).$$
Let now $0<\delta <\pi \cos\tilde{\theta}$ and set
\begin{eqnarray*}
W_{\tilde{\theta}}^F & = & \left\{ w\in \bbC \ \bigg\vert  -\pi \cos\tilde{\theta} + \delta < \begin{pmatrix} \Re(w)\\ \Im(w)\end{pmatrix} \cdot \begin{pmatrix} \cos\tilde{\theta} \\ -\sin\tilde{\theta}\end{pmatrix} < \pi \cos\tilde{\theta} -\delta\right\}\\
& = & \left\{ w\in \bbC \ \bigg\vert  -\pi \cos\tilde{\theta} + \delta < |w| \cos(\tilde{\theta} + \theta_w) < \pi \cos\tilde{\theta} -\delta\right\}
\end{eqnarray*}
and 
$$ U_{\tilde{\theta}}^F  =  \left\{ \gamma \in \bbC \bigg\vert  \tilde{\theta} -\frac{\pi}2 <\Arg(\gamma)  <\frac{\pi}2 +\tilde{\theta}\right\} = 
%\left\{ \gamma \in \bbC \ \bigg\vert  \begin{pmatrix} \Re(\gamma)\\ \Im(\gamma)\end{pmatrix} \cdot \begin{pmatrix} \cos\theta \\ \sin\theta\end{pmatrix}  >0 \right\}\\
 \left\{ \gamma \in \bbC \ \bigg\vert  \cos(\tilde{\theta} - \theta_\gamma)  >0 \right\}.
$$
This definition of $U_{\tilde{\theta}}^F$ guarantees  that the poles of $f^F_\gamma$, which are all on the imaginary axis, when $\gamma$ is positive real, never crosses $\tilde{\Gamma}^{F}_{\tilde{\theta}}$ as the absolute value of the argument of $\gamma$ grows from zero to $\pi$, not including $\pi$. 
\begin{proposition}\label{Fce}
We have for all $t\in \bbR-\{0\}$ and $w\in W^F_{\tilde{\theta}}$ that
$$|K_F(w,e^{i\tilde{\theta}}t)| \leq  \frac{\sqrt{2}}{\sqrt{\pi}(1-e^{-2\pi\cos\tilde{\theta}})} \frac{1}{\delta (1-\alpha)} e^{-\alpha\delta |t|} t^{-2}, \ \ \forall \alpha \in (0,1)$$
and for all $t\in \bbR-\{0\}$ and  $\gamma\in U_{\tilde{\theta}}^F$ that
$$|f_{F,\gamma}(e^{i\tilde{\theta}}t)| \leq \frac{\sqrt{2}}{\sqrt{\pi}(1-e^{-2})} \frac{1}{\cos(\tilde{\theta} - \theta_\gamma)}|\gamma t|^{-1}.$$
\end{proposition}
From this proposition we see that $\delta_2 = 0$ and 
$$ \delta^F = \alpha \delta, \ \ \ c^F_w = \frac{\sqrt{2}}{\sqrt{\pi}(1-e^{-2\pi\cos\tilde{\theta}})}\frac{1}{(1-\alpha) \delta}$$
and
$$ \tilde{c}^F_\gamma = \frac{\sqrt{2}}{\sqrt{\pi}(1-e^{-2})} \frac{1}{\cos(\tilde{\theta} - \theta_\gamma)} \geq \frac{\sqrt{2}}{\sqrt{\pi}(1-e^{-2})} = c^F.$$
\proof
For $t\in \bbR_+$, $\alpha\in (0,1)$  and  $w\in W^F_{\tilde{\theta}}$, using that 
$$(\pi-\Re(w))\cos \tilde{\theta} +  \Im(w)\sin \tilde{\theta} >\delta >0,$$ 
and recalling that $\tilde{\theta} \in (-\frac{\pi}2,\frac{\pi}2)$, we compute 
\begin{eqnarray*} 
|K_F(w,z)| &\leq &2\frac{e^{-t((\pi - \Re(w))\cos\tilde{\theta} + \Im(w)\sin\tilde{\theta})}}{|1-|e^{-2\pi t (\cos\tilde{\theta} + i \sin\tilde{\theta})}||t} \\
&\leq & 2\frac{t  e^{-(1-\alpha) \delta  t} e^{-\alpha \delta  t}}{(1-e^{-2\pi t \cos\tilde{\theta}}) t^2} \\
&\leq & \frac{2e^{-\alpha\delta t} }{\delta (1-\alpha) \sqrt{2\pi}(1-e^{-2\pi \cos\tilde{\theta}})t^2}.
\end{eqnarray*}
Where we have used that
$$\frac{te^{-(1-\alpha)\delta t}}{(1-e^{-2\pi t\cos\tilde{\theta}})}\leq \frac{ 1}{\delta (1-\alpha)\sqrt{2\pi}(1-e^{-2\pi \cos\tilde{\theta}})},$$
which completes the proof for $t\in \bbR_+$. The proof for $t\in \bbR_-$ is completely similar.

For $t\in \bbR_+$ and $\gamma\in U^F_{\tilde{\theta}}$ we have that
\begin{eqnarray*}
 \frac{1}{|\sinh(\gamma e^{i\tilde{\theta}}t)|} & = & 2 \frac{e^{-t |\gamma| \cos(\tilde{\theta} -\theta_\gamma)}}{|1- e^{-2t|\gamma| \cos(\tilde{\theta} -\theta_\gamma)}|}\\
 &  \leq & 2  \frac{t |\gamma| \cos(\tilde{\theta} -\theta_\gamma)e^{-t|\gamma| \cos(\tilde{\theta} -\theta_\gamma) }}{|1- e^{-2t |\gamma| \cos(\tilde{\theta} -\theta_\gamma)}|} \frac{1}{ \cos(\tilde{\theta} -\theta_\gamma)} \frac{1}{|\gamma t|}\\
  & \leq &  \frac{2}{\sqrt{2\pi}(1- e^{-2})} \frac{1}{ \cos(\tilde{\theta} -\theta_\gamma)} \frac{1}{|\gamma t|}.
 \end{eqnarray*}
and the proof for $t\in \bbR_-$ is completely similar.
\eproof

By Proposition \ref{Fce} we conclude for $(w,\gamma)\in W^F_{\tilde{\theta}} \times U^F_{\tilde{\theta}}$ that
$$ g_\gamma^F (w) = \int_{e^{i\tilde{\theta}}(\bbR + i\epsilon)} \frac{e^{wz}}{\sinh(\pi z) z} \frac{1}{\sinh(\gamma z)}dz$$
and
$$ h_m^F (w) = \int_{e^{i\tilde{\theta}}(\bbR + i\epsilon)} \frac{e^{wz}}{\sinh(\pi z) z}  z^m dz$$
are absolutely convergent and thus they are well defined holomorphic functions on $W^F_{\tilde{\theta}}$. Thus we have for $(w,\gamma)\in W_{\tilde{\theta}}^F\times U_{\tilde{\theta}}^F$ that
$$ P_\gamma^{2n}(w) = \frac{1}{4} \sum_{m=0}^n \frac{2 (1 - 2^{2m-1}) B_{2m}}{(2m)!} \gamma^{2m-1}h_{2m-1}^F(w).$$
So we can conclude that for $(w,\gamma)\in W^F_{\tilde{\theta}}\times U^F_{\tilde{\theta}}$ the required estimates (\ref{Mainestimateinfty}) and (\ref{Mainestimateinftyf}) are satisfied along $\tilde{\Gamma}^F_{\tilde{\theta}}$, further we can take $b=1$ in (\ref{cos}) and thus our Theorem \ref{Thm1} therefore applies. We see that in this example
$$ k^F(r) = \frac{1}{\sin(r)}, \ \ \ r\in (0,\pi)$$
and that
$$c'_n = \inf_{r\in(0,\pi)} b_n(r),$$
where
$$ b_n(r) = \frac1{r^{n}} \left( \frac{c^F}{r} + \frac{1}{\sin(r)}\right), \ \ \ r\in (0,\pi).$$
We observe that $b_n(r) \rightarrow \infty$ for $r\rightarrow 0$ and $r\rightarrow \pi$, thus there exist $r_n\in (0,\pi)$ such that
$$ c'_n  = b_n(r_n).$$
Since 
$$ b_n'(r) = - \frac{1}{r^{n+2}} \left( c^F(n+1) + \frac{r n}{\sin(r)} +  \frac{r^2\cos(r)}{\sin^2(r)}\right)$$
we see that
$$  c^F(n+1)\sin^2(r_n) + n r_n \sin(r_n) +  r_n^2\cos(r_n)= 0.$$
We can produce a numerical table of the first few relevant solutions

\begin{table}[h]
\begin{tabular}{|l|l|}
\hline
$2n$ & $r_{2n}$              \\ \hline
$2$  & $2.42067585291066$ \\ \hline
$4$  & $2.64172230058665$ \\ \hline
$6$  & $2.75972744591817$ \\ \hline
$8$  & $2.83308766213237$ \\ \hline
$10$ & $2.88302232511053$ \\ \hline
\end{tabular}
\end{table}

For $n \geq 2$ we have that $b_n'(r)<0$ for  $r \in (0, \frac{3\pi}4)$, hence $r_n\in (\frac{3\pi}{4},\pi)$ for $n \geq 2$ and $r_n \rightarrow \pi$  as $n\rightarrow \infty$.
In particular
$$c'_{2n} \leq \left(\frac{4c^F}{3\pi } + \sqrt{2}\right) \left( \frac23 \right)^{2n}\left( \frac2\pi \right)^{2n}.$$
Further for $n$ positive interger
$$ C^F_{2n} = \frac{4}{c^F\sqrt{(4n-3)\pi}}  \tilde{c}'_{2n},$$
$$ \tilde{c}'_{2n} = \inf_{\alpha \in (0,1)} \frac{\alpha}{1-\alpha} \frac1{\alpha^{2n}} c_{2n}' = \frac{2n}{(1-\frac{1}{2n})^{2n-1}} c_{2n}' .$$
Collecting all constants and applying Theorem \ref{Thm1} we obtain.
\begin{theorem}
We have the following estimates for all positive integers $n$ 
\begin{equation}\label{FE}
\mid \Log(S_\gamma(w))   -   P_\gamma^{2n}(w) \mid \leq \tilde{C}^F_{2n} |\gamma|^{2n} \delta^{-2n}  (2n)! 
\end{equation}
for all $(w,\gamma) \in W_{\tilde{\theta}}^F\times U_{\tilde{\theta}}^F$ and all $\tilde{\theta} \in (-\frac{\pi}{2},\frac{\pi}{2})$ where 
\begin{eqnarray*} \tilde{C}^F_{2n} &=& \frac{4 \sqrt{2} }{\pi \sqrt{4n-3}} \frac{1}{\left(1-\frac{1}{2n}\right)^{2n-1}} \frac{c'_{2n}}{(1- e^{-2\pi \cos\tilde{\theta}})\cos(\tilde{\theta} - \theta_\gamma)}\\
& \leq & \frac{4 \sqrt{2} }{\pi \sqrt{4n-3}} \frac{\sqrt{2} + \frac{4c^F}{3\pi } }{\left(1-\frac{1}{2n}\right)^{2n-1}}  \left( \frac23 \right)^{2n} \frac{1}{(1- e^{-2\pi \cos\tilde{\theta}})\cos(\tilde{\theta} - \theta_\gamma)}\left(\frac2\pi \right)^{2n}\\
&\leq & \frac{3}{(1- e^{-2\pi \cos\tilde{\theta}})\cos(\tilde{\theta} - \theta_\gamma)}\left(\frac2\pi \right)^{2n}.
\end{eqnarray*}
\end{theorem}
The above weakest estimate for $\tilde{C}^F_{2n}$ together with (\ref{FE}) is similar, but not identical, to the estimates obtained, using advanced techniques from the theory of resurgence, in \cite{GK}, in the case where $\Re(\gamma)>0$, e.g. in the case our $\tilde{\theta}=0$. Ours is more general, since it does not require that $\Re(\gamma)>0$, in fact, we see that as we vary $\tilde{\theta} \in (-\frac{\pi}{2},\frac{\pi}{2})$ the above theorem applies to all $\gamma \in \bbC - \bbR_-$. We remark that our general Theorem \ref{Thm1.5} of course also applies and gives a corresponding estimate for  the case $n=0$.

Let us now move on to the applications of Theorem \ref{Thm3} and \ref{Thm5}. 
We recall that for all $z\in\bbC-i\pi\bbZ$
$$\frac{1}{\sinh(z)} =  \frac{1}{z} +  \sum_{n=1}^\infty \left(\frac{(-1)^n }{z - i\pi n} + \frac{(-1)^n }{z + i\pi n} \right)= \frac{1}{z} + 2 \sum_{n=1}^\infty \frac{(-1)^n z}{z^2 + \pi^2 n^2},$$
where the last series converges uniformly on 
$$O_\epsilon = \bbC - \sqcup_{n\in \bbZ-\{0\}} D(i\pi n,\epsilon)$$ 
for all $\epsilon >0$. This means that
$$ \tilde{\phi}^F_{z,p}(\gamma) = \frac{(-1)^n}{\gamma z - i\pi n}$$
and further $m^F=1$, $P_{\phi^F} = i \pi (\bbZ -\{0\})$, $r_m^F = \pi$, $B^F_{i\pi n,1} = (-1)^n$. In relation to condition (\ref{Asump1}) we see that we can take $c=0$ and $\tilde{C} = 2$.
Condition (\ref{Asump2}) is trivially satisfied by
$$\psi^F_z = 0$$
since
$$ \phi^F_z(\gamma) = \frac{1}{\sinh(\gamma z)} - \frac{1}{\gamma z}  =  \frac{2z}{\pi^2} \sum_{n=1}^\infty \frac{(-1)^n}{n^2} \frac{\frac{1}{\gamma}}{(\frac1\gamma)^2 + (\frac{z}{\pi n})^2}$$
converges uniformly in
$$ \gamma \in O_{\epsilon,z} = \bbC - \sqcup_{n\in \bbZ-\{0\}} D(\frac{i\pi n}{z},\epsilon)$$
for $\epsilon>0$ and $z\neq 0$.

For condition (\ref{cosrel}), we first of all constrain $\theta$ and $\gamma$ such that
$$ \cos(\theta - \theta_\gamma) >0.$$
Then we can simply let
$$ r_\theta(w) := \frac{\pi}{|\gamma|}\cos(\theta - \theta_\gamma) =: 2\alpha_\theta(w).$$
If we then require that
$$ \cos(\theta + \tilde{\theta} \pm \frac{\pi}{2}) = 0$$
which is equivalent to
$$ \theta = - \tilde{\theta} + \pi \bbZ,$$
then (\ref{cosrel}) is satisfied, since $\theta_p = \pm \frac{\pi}{2}$ for all $p\in P_{\phi^F}$,
$$ \frac{1}{|\gamma|} \cos(\theta - \theta_\gamma) > \frac{r_\theta(w)}{\pi} > a_\theta(w)$$
and further
$$ \alpha_\theta(w) = \sup_{z\in \overline{\Gamma}^F_w} |z|\cos(\theta + \theta_z \pm \frac{\pi}{2})$$
when $\overline{\Gamma}_w$ agrees with $\tilde{\Gamma}^F_w = \tilde{\Gamma}_{\tilde{\theta}}$ outside the disc $D(0,r_\theta(w)/2)$ and inside the disc $D(0,r_\theta(w)/2)$, it is isotopic  inside $D(0,r_\theta(w)/2) - \{0\}$ (relative to their common boundary) to $\tilde{\Gamma}^F_{\tilde{\theta}}\cap D(0,r_\theta(w)/2)$.
The estimates we developed along $\tilde{\Gamma}^F_w$ of course also holds along $\overline{\Gamma}^{F,\infty}_w(r_\theta(w))$, thus (\ref{cos}), (\ref{Mainestimateinfty}) and (\ref{Mainestimateinftyf}) holds along $\overline{\Gamma}^{F,\infty}_w(r_\theta(w))$. Thus Theorem \ref{Thm5} applies, so we compute
$$ \int_{\overline{\Gamma}_w} K^F(w,z) e^{\frac{z\xi}{i\pi n}} z dz = \int_{\overline{\Gamma}_w} \frac{e^{(w+ \frac{\xi}{i\pi n})z}}{\sinh(\pi z) z} z dz = - 2i \frac{1}{1 + e^{-i(w + \frac{\xi}{i\pi n})}},$$
since, as we observed above, differentiating
$$ \frac{1}{2i} \Li_2(-e^{iw}) = \frac{1}{4} \int_{\Gamma_{w}^{F}}  \frac{e^{wz}}{\sinh(\pi z)  z^2}dz$$
twice in $w$ gives
$$ \frac{1}{1+e^{-iw}} = \frac{i}2 \int_{\Gamma_{w}^{F}} \frac{e^{wz}}{\sinh(\pi z)} dz.$$
Thus Theorem \ref{Thm5} imidiately gives us that
\begin{theorem}  For all  $\theta \in (-\frac{\pi}{2},\frac{\pi}{2})$, $w\in W^F_\theta$ and  $\xi \in V_{w,\theta}$ we have that
$$B_w(\xi)  =  \frac2{i \pi^2} \sum_{n=1}^\infty \frac{(-1)^n}{n^2}  \left(\frac{1}{1+ e^{-i(w+ i\frac{\xi}{\pi n})}} + \frac{1}{1+ e^{-i(w- i\frac{\xi}{\pi n})}}\right)$$
\end{theorem}

We observe that this formula matches the formula obtained in \cite{GK} in the special case where $\theta=0$, e.g. in the case where $\Re(\gamma)>0$.

 From this we can of course now apply Theorem \ref{Thm3} to obtain formula for the Stokes coefficients, since they following immediately from that theorem and the above expression for $B_w$. We observe that
$B_w$ actually extends to a meromorphic function on all of $\bbC$, e.g. $B_w \in \calM(\bbC)$ for all 
$$w \in \bbC - (\pi + 2\pi \bbZ).$$
Furthermore its poles are 
$$ P_{B_w} =  \pm i \pi \bbZ_+(\pi -w + 2\pi \bbZ).$$

\section{Further examples}

\label{Examples}

As further examples of interesting meromorphic transformations we recall a few classical examples, to which our main Theorems applies. %\ref{Thm1}, \ref{Thm1.5}, \ref{Thm2}, \ref{Thm3} and \ref{Thm5}, 
But before we get started on further examples, we make the following observation regarding the Mellin transform and its kernel.
The kernel of the Mellin transform is
$$K_{M}(w,z) = z^w.$$
The generalized Mellin transform with respect to a contour $\Gamma$ is
$$M_{\Gamma}(f)(w) = \int_{\Gamma} K_{M}(w,z) f(z) dz,$$
where we will assume that we can find a holomorphic branch of $K_M$ along $\Gamma$. 
In order to satisfies (\ref{Mainestimateinfty}) we will in some of the examples below consider
$$K_{M,g}(w,z) = z^w g(z)  \ \ \text{ and } \ \ f_g(z) = g(z) f(z).$$
for some $g\in \calM(\bbC)$ such that $P_g \cap \Gamma = \emptyset$ and then of course
$$M_{\Gamma}(f)(w) = \int_{\Gamma} K_{M,g}(w,z) f_g(z) dz.$$
All of the following examples are indeed the Mellin transform (or rescallings thereof) applied to various elementary function, where we have applied above trick for some $g$. 

We will need the Hankel contours
$$ H_-^\epsilon = \{ z\in \bbC \mid d(z, \bbR_-) = \epsilon\} $$
for some $\epsilon >0$, where $H_-^\epsilon$ is orient from $-\infty-i\epsilon$ to $-\infty+ i \epsilon$ and where $d(z,X)$ be the distance from $z$ to $X$. 

Let now $\tilde{H}_-$ be a continuous deformation of $H_-^\epsilon$ inside
$\bbC - \bbR_-$, such that $\tilde{H}_-$ is contained in the half plane $\{z \in \bbC  | \Re(z) < 0\}$ except it passes through $0$. Further we will assume that (\ref{cos}), with $b=\frac12$, is satisfies along the two smooth pieces of
$\tilde{H}_-$ inside $\{z \in \bbC | \Re(z) < 0\}$.

\subsection{Euler's Gamma function $\Gamma$}

We now let 
$$ K^\alpha_\Gamma(w,z) = \frac{iz^{w-1}e^{\alpha z }}{2\sin(\pi w)} \text{ and } f_{\Gamma,\alpha}(z) = e^{(1-\alpha) z }$$
for $\alpha \in (0,1)$ and 
$$ \Gamma^\Gamma_{w} = H_-^\epsilon. $$
We then have that
$$ \Gamma^\alpha_\gamma(w)=  g^{\Gamma,\alpha}_\gamma(w)$$
where $\Gamma^\alpha_1$ is Euler's Gamma function for all $\alpha \in (0,1)$.

We now observe that for any $0<\delta_1 < \alpha$ we that (\ref{Mainestimateinfty}) is satisfied along
$$ \tilde{\Gamma}^\Gamma_{w} = \tilde{H}_-$$
and that for all $\delta_2 > 1-\alpha$ we have that (\ref{Mainestimateinftyf}). But since we need $\delta_2 < \delta_1$ we see that we need $\alpha \in (\frac12,1)$. Then Theorem \ref{Thm1} gives an asymptotic expansion
\begin{equation}\label{gamma1}
 g^{\Gamma,\alpha}_\gamma(w) \sim \sum_{m=0}^\infty h^{\Gamma,\alpha}_m(w) \frac{(1-\alpha)^m}{m!} \gamma^m
\end{equation}
where
$$ h^{\Gamma,\alpha}_m(w) = \frac{\partial^m}{\partial \alpha^m} h^{\Gamma,\alpha}_0(w)$$
and 
$$h^{\Gamma,\alpha}_0(w) = \frac{i}{2\sin(\pi w)} \int_{H^\epsilon_-} e^{\alpha z} z^{w-1} dz.$$
By Proposition \ref{PropConv}, we actually have that the right hand of (\ref{gamma1}) is convergent for $|\gamma| <\alpha R$ and equal to the left hand side, where
$$ R = \frac{1}{1-\alpha},$$
since
$$ \left( \frac{(1-\alpha)^m}{m!} m!\right)^{1/m} = 1-\alpha.$$

\subsection{The reciprocal of Euler's Gamma function $\frac{1}{\Gamma}$}

We now let 
$$ K^\alpha_{\frac{1}{\Gamma}}(w,z) = \frac{1}{2\pi i}z^{-w} e^{\alpha z} \text{ and } f_{\frac{1}{\Gamma},\alpha}(z) = e^{(1-\alpha)z}$$
for $\alpha \in (0,1)$ and 
$$ \Gamma^{\frac{1}{\Gamma}}_{w} = H^\epsilon_-.$$
We then have that
$$ \frac{1}{\Gamma^\alpha_\gamma(w)} = g^{\frac{1}{\Gamma},\alpha}_\gamma(w)$$
for $\alpha \in (0,1)$.
Again we can apply Theorem \ref{Thm1} for $\alpha \in (\frac12,1)$ to get the asymptotic expansion
\begin{equation}\label{gamma2}
 g^{\frac{1}{\Gamma},\alpha}_\gamma(w) \sim \sum_{m=0}^\infty h^{\frac{1}{\Gamma},\alpha}_m(w) \frac{(1-\alpha)^m}{m!} \gamma^m
\end{equation}
where
$$ h^{\frac{1}{\Gamma},\alpha}_m(w) = \frac{\partial^m}{\partial \alpha^m} h^{\frac{1}{\Gamma},\alpha}_0(w)$$
and 
$$h^{\frac{1}{\Gamma},\alpha}_0(w) = \frac{1}{2\pi i} \int_{H^\epsilon_-} e^{\alpha z} z^{-w} dz.$$
By Proposition \ref{PropConv}, we actually have that the right hand of (\ref{gamma2}) is convergent for $|\gamma| <\frac{\alpha}{1-\alpha} $.

\subsection{Riemann zeta function $\zeta$}

We set 
$$ K_\zeta(w,z) = -\frac{\Gamma(1-w)}{2\pi i}  z^{w-1} e^{\alpha{z}}\text{ and } f_\zeta(z) = \frac{e^{(1-\alpha){z} }}{1- e^{z}}$$
for $\alpha \in (0,1)$ and 
$$ \Gamma^\zeta_{w} = H^\epsilon_-.$$
We then have that for $w\neq 1$
$$ \zeta^\alpha_\gamma(w) = g^{\zeta,\alpha}_\gamma(w)$$
where $\zeta = \zeta^\alpha_1$ for $\alpha \in (0,1)$.

Again we can apply Theorem \ref{Thm1} for $\alpha \in (\frac12,1)$ to get the asymptotic expansion
\begin{equation}\label{zeta1}
 g^{\zeta ,\alpha}_\gamma(w) \sim \sum_{m=-1}^\infty h^{\zeta ,\alpha}_m(w) a_m \gamma^m
\end{equation}
where
$$ h^{\zeta,\alpha}_m(w) = \frac{\partial^m}{\partial \alpha^m} h^{\zeta ,\alpha}_0(w), \ \ \ a_m = \sum_{m_1=0}^m \frac{(1-\alpha)^{m_1}}{m_1!} \frac{B_{m-m_1+1}}{(m-m_1+1)!}$$
and 
$$h^{\zeta ,\alpha}_0(w) = -\frac{\Gamma(1-w)}{2\pi i} \int_{H^\epsilon_-}  z^{w-1} e^{\alpha z}dz = \alpha^{-w}.$$
We leave it to the reader to check the applicability of Theorem \ref{Thm2}, \ref{Thm3} and \ref{Thm5} in this case.

\subsection{Hurwitz zeta function $\zeta$}

We set 
$$ K_H(w,z) = \frac{\Gamma(1-w)}{2\pi i} \frac{e^{(w -1)\log(z)} }{ e^{-z}-1}\text{ and } f_H(z) = e^{z}$$
and 
$$ \Gamma^H_{w} = H^\epsilon_-.$$ 
We observe that the Hurwitz zeta functions equals
$$ \zeta^H(w,q) = g^H_{q-1}(w).$$

Let now $\tilde{\Gamma}^H_w = \tilde{H}_-$. We compute that
$$ \bigg\vert \frac{e^{(w -1)\log(z)} }{ e^{-z}-1} \bigg\vert = \frac{e^{(\Re(w)-1)\log|z| - \Im(w) \text{Arg}(z)} e^{\Re(z)}}{|1-e^z |}$$
Thus there exist positive constants $C_{\delta_1}$ for all $\delta_1<1$ and $c_w$ for all $w\in \bbC$, such that
$$ \vert K_H(w,z) \vert  \leq C_{\delta_1} c_w e^{-\delta_1|z|}|z|^{-1}$$
for all $z\in \tilde{\Gamma}_w^H$ and all $w \in \bbC$. Further we see that
$$ \tilde{\phi}_n(\gamma z) = \sum_{k=n+1}^\infty \frac{|\gamma z|^{k-n-1}}{k!} \leq e ^{|\gamma| |z|}  \leq e^{\tilde{\delta} |z|}$$
also for all $z\in \tilde{\Gamma}_w^H$ whenever $|\gamma| \leq \tilde{\delta}$. Thus Theorem \ref{Thm77} applies to establish that $g_\gamma^H(w)$ has an asymptotic expansion
$$  g_\gamma^H(w) \sim \sum_{m=0}^\infty \zeta_m(w) \frac{\gamma^m}{m!}$$
near $\gamma = 0$ where
$$ \zeta_m(w) = \int_{\Gamma^H_w} K_H(w,z) z^m dz.$$ 
In fact Theorem \ref{Thm77} provides the estimate
$$  \vert g_\gamma^H(w) - \sum_{m=0}^n \zeta_m(w) \frac{\gamma^m}{m!} \vert \leq 4 C_{\delta_1} c_w |\gamma|^{n+1} \delta^{n+1} n!$$
for $\gamma \in U_{\tilde{\delta}} = D(0,\tilde{\delta})$, where
$$ \delta= \delta_1 - \tilde{\delta} >0.$$
By proposition Proposition \ref{PropConv} and since
$$ \left(\frac{(m-1)!}{m!}\right)^{\frac1m} \rightarrow 1 \ \ \text{ for } \ \ m \rightarrow \infty$$ we see in fact that the asymptotic expansion is convergent for $|\gamma| < \delta_1$ and we thus have that
$$  g_\gamma^H(w) = \sum_{m=0}^\infty \zeta_m(w) \frac{\gamma^m}{m!}$$
for $|\gamma| < \delta_1$. We observe that 
$$ \zeta(w) = \zeta_0(w).$$

\subsection{Gauss hypergeometric function $ _2 F_1$} 
We set
$$ K_{_2 F_1, \alpha}(w,z) =  z \frac{e^{z^2 \log(-w)}}{\cos(\alpha z^2)} \text{ and } f_{_2 F_1,\alpha}(z) = \cos(\alpha z^2) \frac{\Gamma(-z^2) \Gamma(z^2+a) \Gamma(z^2+b)}{\Gamma(z^2+c)} $$
 for $a,b$ non-positive integers and
$$\Gamma^{_2 F_1,\epsilon} = \left((e^{i\frac{3\pi}{4}} \bbR_- \cup e^{i\frac{\pi}{4}} \bbR_+)\cap (\bbC - D(0,\epsilon))\right) \cup \{ z \mid |z| = \epsilon, -\frac{\pi}{4} < \text{Arg}(z) < \frac{\pi}{4}\} $$
with the orientation induced from the usual one on $\bbR_\pm$ and the positive orientation on the circle. Then
$$ \frac{\Gamma(a)\Gamma(b)}{\Gamma(c)}  {_2 F^{\gamma,\alpha}_1} (a,b;c;w) = \frac{1}{\pi i} g_\gamma^{_2 F_1,\alpha}(w)$$
with $ _2 F^{1,\alpha}_1(a,b;c;w) =  {_2 F_1}(a,b;c;w) $ Gauss hypergeometric function for all $\alpha \in (0,\pi)$.
Let 
$$ \tilde{\Gamma}^{_2 F_1,\epsilon} = e^{i\frac{3\pi}{4}} \bbR_- \cup e^{i\frac{\pi}{4}} \bbR_+ .$$
We observe that for any $\tilde{\delta} <\pi$ there exist $C_{\tilde{\delta},\alpha}$ such that
$$ |f_{_2 F_1,\alpha}(z) | \leq C_{\tilde{\delta},\alpha} e^{-(\tilde{\delta} - \alpha) |z|^2}|z|^{-2} \ \ \ \forall z\in \tilde{\Gamma}^{_2 F_1}_w$$
and
$$ |K_{_2 F_1, \alpha}(w,z)| \leq C_{\tilde{\delta},\alpha} e^{-(\alpha - |\text{Arg}(w)|) |z|^2} \ \ \ \forall z\in \tilde{\Gamma}^{_2 F_1}_w$$
Thus if we assume that
$$ \tilde{\delta} >\alpha \ \ \ \text{ and } \ \ \ |\text{Arg}(w)|  < \alpha $$
we can let $\delta_1 = \tilde{\delta} - \alpha$ and $\delta_2=0$ and then we satisfy (\ref{Mainestimateinfty}) and (\ref{Mainestimateinftyf}). It is clear that (\ref{cos}) is satisfies with $b=1$. Thus Theorem \ref{Thm1} applies to establish that we get an asymptotic expansion
\begin{equation}\label{GH1}
 g^{_2 F_1,\alpha}_\gamma(w) \sim \sum_{m=-2}^\infty h^{_2 F_1,\alpha}_m(w) a_m \gamma^{m}
\end{equation}
where $a_m$ are the coefficients of the Laurent expansion of $ f_{_2 F_1,\alpha}$ around $0$ and 
$$h^{_2 F_1,\alpha}_m(w) =  \int_{\Gamma^{_2 F_1,\epsilon}} \frac{e^{z^2 \log(-w)}}{\cos(\alpha z^2)} z^{m+1} dz$$
for $m=-2,-1,0,1,2, \ldots.$ We leave it to the reader to check the applicability of Theorem \ref{Thm2}, \ref{Thm3} and \ref{Thm5} in this case.

\subsection{The Airy function $\text{Ai}$} 
We set
$$ K_{\text{Ai},a}(w,z) = e^{-zw} e^{\alpha z} \text{ and }  f_{\text{Ai}}(z) = e^{\frac{1}{3} z^3 - \alpha z}$$
for $a\in (0,\frac13)$ and
$$\Gamma^{\text{Ai}} = \bbR^{\pi-\tilde{\theta}}_- \cup \bbR_+^{\tilde{\theta}}$$
where $\tilde{\theta} \in [\frac{\pi}{6}, \frac{\pi}{3}]$ and $\bbR^{\tilde{\theta}}_- = \bbR_- e^{i\tilde{\theta}}$.
This implies  that
$$ \cos(3\tilde{\theta}) < 0$$ 
thus (\ref{Mainestimateinfty}) for any $\delta_1 \in \bbR_+$ such that
$$ |w| \cos(\tilde{\theta} + \theta_w) - \alpha \cos(\tilde{\theta}) > \delta_1$$
and (\ref{Mainestimateinftyf}) with $\delta_2 = 0$ are satisfied, hence Theorem \ref{Thm1} applies since (\ref{cos}) is obviously satisfied with $b=1$.
Thus we get an asymptotic expansion
\begin{equation}\label{Airy1}
 g^{\text{Ai},\alpha}_\gamma(w) \sim \sum_{m=0}^\infty h^{\text{Ai},\alpha}_m(w) a_m \gamma^{m}
\end{equation}
where
$$ h^{\text{Ai},\alpha} _m(w) = \frac{\partial^m}{\partial \alpha^m} h^{\text{Ai},\alpha}_0(w), \ \ \  \ \ \ a_m = \sum_{l=0}^{\frac{m}3} (-1)^{m-l} \frac{\alpha^{m-3l}}{l! (m-3l)! 3^l}$$
and 
$$h^{\text{Ai} ,\alpha}_0(w) =  \int_{\Gamma^{\text{Ai}}} e^{-zw} e^{\alpha z} dz.$$
By Proposition \ref{PropConv}, we actually have that the right hand of (\ref{Airy1}) is convergent for all $| \gamma | <\delta_1 R$, where
$$ R\geq  \frac1{1 +\alpha}.$$

\end{document}